\documentclass[11pt]{article}

\usepackage{fullpage,amsfonts,amsmath,amsthm,amssymb,mathrsfs,graphicx,epstopdf}
\usepackage[top=2.54cm, bottom=2.54cm, left=2.54 cm, right=2.54 cm]{geometry}

 \newcommand{\lab}[1]{\label{#1}}                

 \usepackage[usenames,dvipsnames]{color}

\newcommand{\remove}[1]{}
\newcommand\eqn[1]{(\ref{#1})}

\newcommand{\be}{\begin{equation}}
\newcommand{\bel}[1]{\begin{equation}\lab{#1}\ }
\newcommand{\ee}{\end{equation}}
\newcommand{\bea}{\begin{eqnarray}}
\newcommand{\eea}{\end{eqnarray}}
\newcommand{\bean}{\begin{eqnarray*}}
\newcommand{\eean}{\end{eqnarray*}}

\newtheorem{thm}{Theorem}
\newtheorem{cor}[thm]{Corollary}

\newtheorem{lemma}[thm]{Lemma}
\newtheorem{definition}[thm]{Definition}
\newtheorem{claim}[thm]{Claim}

\newtheorem{remark}[thm]{Remark}
\newtheorem{example}[thm]{Example}

\def\proof{\noindent{\bf Proof.~}~}
\def\qed{~~\vrule height8pt width4pt depth0pt}

\newcommand{\norm}[1]{{\left\|#1\right\|}}

\newcommand{\ind}[1]{1_{\{#1\}}}
\def\var{{\bf Var}}

\def\tbfd{\widetilde{\bfd}}
\def\td{\widetilde{d}}
\def\tM{\widetilde{M}}

\def\B{{\mathcal B}}
\def\C{{\mathcal C}}

\def\F{{\mathcal F}}
\def\G{{\mathcal G}}

\def\I{{\mathcal I}}
\def\J{{\mathcal J}}

\def\N{{\mathcal N}}

\def\W{{\mathcal W}}

\def\ex{{\mathbb E}}
\def\pr{{\mathbb P}}

\def\bfd{{\bf d}}

\def\bfg{{\bf g}}
\def\bfh{{\bf h}}
\def\bfi{{\bf i}}
\def\bfj{{\bf j}}

\def\eps{\epsilon}

\newcommand{\fix}[1]{{#1}}

\date{}

\title{Triangles and subgraph probabilities in random regular graphs}
\author{Pu Gao\thanks{Research supported by NSERC.} \\
University of Waterloo\\
pu.gao@uwaterloo.ca }
\begin{document}
\maketitle

\begin{abstract}

We improve the estimates of the subgraph probabilities in a random regular graph. Using the improved results, we further improve the limiting distribution of the number of triangles in random regular graphs.

\end{abstract}

\section{Introduction}

Research in random graph theory started from the study of subgraphs~\cite{erdHos1960}. The distributional results for small subgraphs and for large subgraphs are very different in nature. 
The distributions of small subgraphs in $\G(n,p)$ are well understood~\cite{ruciniski1985,rucinski1988}.  Under some mild conditions, the number of subgraphs isomorphic to a fixed graph $H$, denoted by $Z_H$, is  asymptotically normally distributed when $p$ exceeds some critical value $p_H$. Similar results hold for $\G(n,m)$ as well~\cite{ruciniski1985}. In particular, if $H$ is a balanced graph and $p$ is not too close to 1, then $Z_H$  is asymptotically normally distributed in $\G(n,p)$, and is thus highly concentrated around its expectation $\ex Z_H$, when $\ex Z_H\to\infty$ as $n\to\infty$.   The picture for large subgraphs is very different. Take $Z_M$, the number of perfect matchings (assuming $n$ is even) for an example. The threshold of its appearance is at $p=\log n/n$. However, $\ex Z_M\to\infty$ already when $p=C/n$ for $C>e$. This implies no concentration of   $Z_M$ around its expectation for $p=O(1/n)$, and indeed no such concentration for larger $p$ either. It has been proved that $Z_M$ is asymptotically log-normally distributed for $p\gg n^{-1/2}$~\cite{janson1994numbers}. For smaller values of $p$ its distribution is unknown.  Similar results have been proved for the number of spanning trees and the number of Hamilton cycles~\cite{janson1994numbers}, the number of $d$-factors~\cite{gao2013distribution}, and the number of triangle factors~\cite{gao2013distributions}. What happens to the number of subgraphs whose size is between constant size and size of order $\Omega(n)$? For the number of matchings of size $\ell$, it turns out that its limiting distribution transitions from normal to log-normal at some critical point of $p$~\cite{gao2016transition}, and this distributional transition may be a general phenomenon for other subgraphs as well. 

Another extensively studied random graph model is the model of random regular graphs. Let $\G(n,d)$ denote a graph chosen uniformly at random from the set of all $d$-regular graph on vertex set $[n]$. Much less is known about the distribution of $Z_H$ in $\G(n,d)$ where $H$ is of fixed size, or of size $\Omega(n)$. For $d=O(1)$, $\G(n,d)$ is locally tree like and thus, the only nontrivial  connected subgraphs of fixed sizes are cycles (see~\cite{wormaldsurvey} for more discussions of properties of $\G(n,d)$ for fixed $d$).  
In $\G(n,d)$ the joint distribution of the number of cycles of any fixed length tends asymptotically to the distribution of  independent Poisson random variables~\cite{bollobas1980probabilistic,wormald1981asymptotic}.  This continues to hold for slowly growing $d$~\cite{mckay2004short}. For instance, the Poisson like distribution for the number of triangles is proved to hold for all $d=o(n^{1/5})$. Z. Gao and Wormald~\cite{gao2008distribution} proved normal distribution for strictly balanced graphs, which, in cases of short cycles, permits slightly larger $d$ than in~\cite{mckay2004short}. For the number of triangles, their result holds for $d=o(n^{2/7})$. For large subgraphs, Janson~\cite{janson1993random} determined the limiting distribution of the number of perfect matchings and the number of Hamilton cycles in $\G(n,d)$ for constant $d$. For growing $d$, we are not aware of any results when the size of $H$ is beyond $\log n$. The recent development on the sandwich conjecture~\cite{gao2020sandwiching,gao2020Kim} allows to translate properties and graph parameters from $\G(n,p)$ to $\G(n,d)$, but in general the distributional results do not translate.

There are two main obstacles in proving the limiting distribution of $Z_H$, even for $H$ of fixed size. The first obstacle lies in the computation of the moments of $Z_H$. To prove normality, one can compute its central $j$-th moment for each fixed $j$. However, due to the appearance of both positive and negative terms in the central moments computation, the leading terms cancel each other, and the errors in the subgraph probabilities quickly become the dominating terms when $d$ grows. (Take the triangle $H=C_3$ for an example. We  only manage to compute the variance of $Z_{C_3}$ for $d=o(n^{2/5})$, and computing higher moments requires even smaller $d$.) The other standard method is computing the $k$-th raw moment for $k$ around $\ex Z_H/\sigma (Z_H)$, where $\sigma (Z_H)$ is the standard deviation of $Z_H$. For such large $k$ (which typically grows quickly with $d$), it is difficult to bound the number of $k$-tuples of copies of $H$ sharing a specific number of edges. The second obstacle is to  estimate the probabilities of subgraphs in $\G(n,d)$. Specifying a set of edges $H\subseteq \binom{[n]}{2}$, what is the probability that $H\subseteq \G(n,d)$ (i.e.\ $\G(n,d)$ contains all edges in $H$)? McKay~\cite{mckay1981subgraphs,mckay2010subgraphs} gave an estimate of this probability which has a relative error $O(|H|d/n)$, and to our knowledge, this has been the best estimate to date. This immediately limits the study of $Z_H$ for large $H$, or the computation of high moments of $Z_H$ even for $H$ of constant size, which is required when proving the limiting distribution of $Z_H$. Take $H=C_{3}$ as an example. To determine the limiting distribution of $Z_{C_3}$, we need to estimate the $k$-th moment of $Z_{C_3}$ where $k\approx d^{3/2}$. This requires an asymptotic joint probability estimate of $\Theta(d^{3/2})$ edges. With errors in~\cite{mckay1981subgraphs,mckay2010subgraphs}, $d$ is required to be $o(n^{2/5})$.  

For the  second obstacle, we manage to improve the error in the subgraph probabilities in~\cite{mckay1981subgraphs,mckay2010subgraphs}. See Theorem~\ref{thm:conditional} in Section~\ref{sec:results} for the precise statement. This new result allows us to estimate, for instance, the asymptotic probability that $\G(n,d)$ contains a specific perfect matching for $d$ up to $n^{1/2}$. As we only study small subgraphs in this paper, it turns out that a less precise form would be sufficient, which is stated in Theorem~\ref{thm:special} as a corollary of Theorem~\ref{thm:conditional}.

For the first obstacle, Z. Gao and Wormald bypassed it by smartly considering a new variable ${\widetilde Z}_H$, which counts isolated copies of $H$, i.e.\ copies of $H$ that do not share edges with other copies of $H$. In the range of $d$ where the number of non-isolated copies of $H$ is so small that it does not affect the standard deviation of ${\widetilde Z}_H$, the distribution of ${\widetilde Z}_H$ will immediately yield the distribution of $Z_H$. Working on ${\widetilde Z}_H$ avoids having to deal with tuples of heavily intersecting copies of $H$. However, it also limits the range of $d$ for which the proof method can apply. For the case $H=C_3$, $d=o(n^{2/7})$ is required in~\cite{gao2008distribution}, which did not even reach the natural stopping point $d=o(n^{2/5})$ (explained above). To relax the condition on $d$, it is possible to modify the definition of ${\widetilde Z}_H$. Instead of counting isolated triangles,  we may, for instance, count triangles that intersect only a bounded number of other triangles. With the new definition and with the new subgraph probabilities in   Theorem~\ref{thm:special}, it is possible to get beyond $d=o(n^{2/7})$. However it is not sufficient to reach $d=n^{1/2-o(1)}$. In this paper, we handle $Z_H$ directly, and prove the normality of the number of triangles for all $d=O(n^{1/2})$ where $d\to\infty$.

We also computed the variance of $Z_H$ when $H$ is a fixed strictly balanced graph. Surprisingly, although the variance seemingly carries less information than the limiting distribution and thus computing it should be easier than determining the limiting distribution, we are not able to determine the variance of the number of triangles for $d$ up to $n^{1/2}$. In fact, we can only compute the variance for $d=o(n^{2/5})$. Note that the limiting distribution does not imply the variance, as it is not sensitive to a lottery effect (events affecting $Z_H$ that occur with a tiny probability), whereas the variance is, although we do not expect such  a lottery effect for the variance either.  Another surprising discovery is that the number of triangles seems to behave more like the sum of a set of independent variables in $\G(n,d)$ than in $\G(n,p)$, whereas the intuition is the other way around. Let $X_{ijk}$ be the indicator variable that the three vertices $i,j,k$ induce a triangle. Then the number of triangles is $\sum X_{ijk}$ where the sum is over all $3$-subsets of $[n]$. In $\G(n,p)$,  $X_{ijk}$ is independent of all but those $X_{uvw}$ such that $|\{i,j,k\}\cap \{u,v,w\}|=2$. In $\G(n,d)$, due to the dependency between edges, all these indicator variables are correlated. However, surprising cancellations of leading terms appear in the calculations for the variance of the number of triangles (see Section~\ref{sec:variance-triangle} for details). The variance of the number of triangles turns out smaller in $\G(n,d)$ than in $\G(n,d/n)$. 

One of the challenges in computing the moments of the number of triangles is to bound the number of $k$-tuples of triangles where $k\approx d^{3/2}$, such that these $k$ triangles induce a large number of ``holes'' (see Definition~\ref{def:hole} in Section~\ref{sec:small}). The bound and the arguments in this paper do not easily extend to other strictly balanced graphs in general. Extending the distributional result for the number of triangles to other strictly balanced graphs requires a smarter treatment for the holes.

\section{Main results}
\label{sec:results}

\subsection{Improved subgraph probabilities}

 Let $\bfd=(d_1,\ldots,d_n)$ be a degree sequence and let $\G(n,\bfd)$ denote a uniformly random graph with degree sequence $\bfd$. Let $H_1$ and $H_2$ be two disjoint graphs on $[n]$, i.e.\ $E(H_1)\cap E(H_2)=\emptyset$. Let $H_1^+$ denote the event that $H_1\subseteq \G(n,\bfd)$ and let $H_2^-$ denote the event that $H_2\cap \G(n,\bfd)=\emptyset$. The conditional edge probability $\pr(uv\in\G(n,\bfd)\mid H_1^+,H_2^-)$ was first given by  McKay~\cite{mckay1981subgraphs,mckay2010subgraphs}, which immediately applies to give an asymptotic estimate of $\pr(H\subseteq \G(n,\bfd))$ when $H$ is not too large. 
 This result was recently extended by Ohapkin and the author~\cite{gao2020subgraph} to  more general degree sequences, with the same order of the error term $O(\Delta(\bfd)/n)$ as in~\cite{gao2020subgraph} ($\Delta(\bfd)$ denotes the maximum degree in $\bfd$). Given a graph $H$ on $[n]$, let  $\bfd^{H}=(d_1^{H},\ldots, d_n^{H})$ denote the degree sequence of $H$ and let $|H|$ denote the number of edges in $H$. Let $\Delta_H$ denote the maximum degree of $H$. Given two degree sequences $\bfd$ and $\bfd'$, we say $\bfd\preceq \bfd'$ if $d_i\le d'_i$ for every $1\le i\le n$.
Let $\Delta=\Delta(\bfd)$ be the maximum component of $\bfd$  and let $M=\sum_{i=1}^n d_i$. The following result follows as a corollary of~\cite[Theorem 1]{gao2020subgraph}.
 
\begin{thm}\label{thm0:conditional} Let $\bfd$ be such that $\Delta^2=o(M)$.
Let $H_1$ and $H_2$ be two disjoint graphs on $[n]$ where $\bfd^{H_1}\preceq \bfd$. Suppose that $\Delta_{H_2}=O(\Delta)$ and $M-2|H_1|=\Omega(M)$. Suppose further that
 $uv\notin H_1\cup H_2$.  
Then,
\begin{eqnarray*}
\pr(uv\in \G(n,\bfd) \mid H_1^+, H_2^-)& =& \left(1+O\left(\frac{\Delta^2}{M}\right)\right)\frac{(d-d^{H_1}_u)(d-d^{H_1}_v)}{M-2|H_1|}.
\end{eqnarray*}\
\end{thm}

\fix{One of the main results of this paper is an estimate of the above conditional probability with a relative error  $\Delta^5n/M^3$ (See Theorem~\ref{thm:conditional}).  Due to the technical statement of the theorem, we start our discussion from a special and simpler case where $\bfd=(d,\ldots, d)$ and $H_2=\emptyset$.} 
Given graph $H$ and $u\in [n]$, let $\N_H(u)$ denote the set of vertices that are adjacent to $u$ in $H$. The theorem below gives an approximation of $\pr(uv\in \G(n,d)\mid H^+)$ with a less sharp error compared to Theorem~\ref{thm:conditional} below. However, it is sufficient for  studying $Z_H$ in $\G(n,d)$ when $H$ is not too big.

\begin{figure}
  \[
  \includegraphics[scale=0.8]{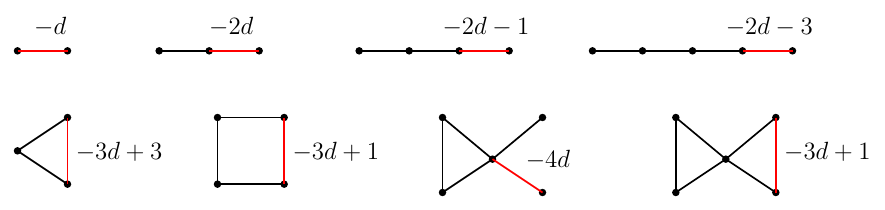}  
  \]
  \caption{$\phi_{H}(uv)$}
  \label{f:phi}
\end{figure}

\begin{thm}\label{thm:special}
Let $H$ be a graph on $[n]$ where $d_i^{H}\le d$ for every $i$. Suppose that $d=o(n)$, $dn-|H|=\Omega(dn)$, and that
 $uv\notin H$. Let $\tbfd=\bfd-\bfd^H$ where $\bfd=(d,\ldots,d)$. Then,
\begin{equation}
\pr(uv\in \G(n,d)\mid H^+)=\frac{\td_u\td_v}{dn}\left(1-\frac{\phi_H(uv)}{dn}\right)\left(1+O\left(\fix{\frac{|H|}{n^2}+\frac{|H|^2}{d^2n^2}+\frac{d^2}{n^2}}\right)\right),\label{phi0}
\end{equation}
where 
\begin{eqnarray}
\phi_H(uv)&=&-d-2|H|-(d-1)(d_u^H+d_v^H)+d_u^Hd_v^H+\sum_{x\in {\cal N}_H(u)}  d^H_x+ \sum_{y\in  {\cal N}_H(v)}d^H_y.\label{phi1}
\end{eqnarray}
\end{thm}
\begin{remark} By the definition of $\phi_H(uv)$ it is immediate that \fix{if the maximum degree of $H$ is $O(1)$ then}
\begin{align}
\phi_H(uv)
&=-2|H|-d(1+d_u^H+d_v^H)+O(1),\label{phi2}
\end{align}
and this bound is  sufficiently accurate in many applications. Given $H$, $\phi_H(uv)$ can be easily computed by examining the neighbourhood of $u$ and $v$ in $H$.
In Figure~\ref{f:phi} we give a few examples of the value of $\phi_H(uv)$. Edge $uv$ is coloured red, and edges in $H$ are coloured black. The value of $\phi_H(uv)$ is given next to the edge $uv$.
\end{remark}

\begin{remark} \fix{We check that $\pr(uv\in \G(n,d))$ from Theorem~\ref{thm:special} agrees with what we know about this probability.}
By Theorem~\ref{thm:special}, 
\[
\pr(uv\in \G(n,d))=\fix{\frac{d}{n} \left(1+\frac{d}{dn}\right) } (1+O(d^2/n^2))=\frac{d}{n-1}(1+O(d^2/n^2)),
\] 
whereas by symmetry, we know that $\pr(uv\in \G(n,d))=d/(n-1)$. Hence the estimate from the theorem is accurate up to a relative error $O(d^2/n^2)$, as it is supposed to.
\end{remark}

The following corollary on the upper bound of $\pr\big(uv\in\G(n,d)\mid H^+\big)$ will be useful in the study of the distribution of the number of triangles. The proof of the corollary is given in Section~\ref{sec:corollary}.

 \begin{cor}\label{cor:conditional}
Let $\bfd=(d,d,\ldots, d)$ where $dn$ is even and $d=o(n)$.
Let $H$ be a graph on $[n]$ with degree sequence $\bfd^H\preceq \bfd$ such that $dn-2|H|=\Omega(dn)$. 
Suppose $F\subseteq K_n\setminus H$ where $|F|=O(1)$. Let $\lambda_F=\pr(F\subseteq \G(n,d))$ and $\tbfd=\bfd-\bfd^{H}$. Then,
\begin{eqnarray}
\pr\Big(F\subseteq\G(n,d)\mid H^+\Big)
&\le& \lambda_F\left(1+O\left(\frac{1}{n}+\frac{|H|}{dn}+\frac{d^2}{n^2}\right)\right).
\end{eqnarray}
\end{cor}

Theorem~\ref{thm:special} is a special and less precise version of the following more general result. Given $\bfd$, recall that $M=M(\bfd)=\sum_{i=1}^n d_i$. Let $M_j=M_j(\bfd)=\sum_{i=1}^n (d_i)_j$ for any integer $j\ge 2$.

\begin{thm}\label{thm:conditional} Let $\bfd=(d_1,\ldots,d_n)$ be such that $M$ is even and $\Delta^2=o(M)$, where $\Delta=\Delta(\bfd)$. 
Let $H_1$ and $H_2$ be disjoint graphs with degree sequences $\bfd^{H_1}\preceq \bfd$ and $\Delta_{H_2}=O(\Delta)$. Let $\tbfd=\bfd-\bfd^{H_1}$ and let $\tM=M(\tbfd)$ and $\tM_j=M_j(\tbfd)$ for $j\ge 2$. Suppose $\tM=\Omega(M)$, and $uv\notin H_1\cup H_2$. 
Then
\begin{align*}
&\pr\Big(uv\in\G(n,\bfd)\mid H_1^+, H_2^-\Big)\\
&=\left(1+O\left(\fix{\frac{\Delta^5n}{M^3}}\right)\right)\frac{\td_u\td_v}{\tM}\left(\fix{1-\frac{\tM_2^2}{\tM^3}-\frac{\tM_2}{\tM^2}-\sum_{(x,y)\in \W} \frac{\td_x\td_y}{\tM^2}} \right)\left(1-\frac{\bar\phi_{H_1,H_2}(uv)}{\tM}\right),
\end{align*}
where
\[
\bar\phi_{H_1,H_2}(uv)=-2(\td_u+\td_v)+2-\sum_{x\in {\cal N}_{H_1\cup H_2}(u)} \td_x-\sum_{y\in {\cal N}_{H_1\cup H_2}(v)} \td_y-\frac{(\td_u+\td_v-2)\tM_2}{\tM} +\td_u\td_v ,\label{barphi}
\]
and $\W$ is the set of $(x,y)$ such that $xy\in H_1\cup H_2$, $xu,yv\notin H_1\cup H_2$, and $\{x,y\}\cap\{u,v\}=\emptyset$. 

\end{thm}

\begin{remark}
To study $Z_H$ for large $H$, for instance,  if $H$ is a perfect matching, then the error $O(|H|^2/d^2n^2)$ in Theorem~\ref{thm:special} is too large, as the cumulative error becomes $O(|H|^3/d^2n^2)=O(n/d^2)$ which is $\Omega(1)$ for $d=O(\sqrt{n})$. In such cases Theorem~\ref{thm:conditional} needs to be applied. 
\end{remark}

\fix{
\begin{remark}
A result by McKay~\cite[Theorem 4.6]{mckay1985} applies to find the probability $\pr(H\subseteq \G(n,d))$ with a relative error $d^3/n$. Theorem~\ref{thm:conditional} applies to find this probability with a smaller relative error $d^2|H|/n^2$.
\end{remark}
}

\subsection{Proof of Theorem~\ref{thm:special}}
We prove that Theorem~\ref{thm:special} follows as a corollary of Theorem~\ref{thm:conditional}.
Apply Theorem~\ref{thm:conditional} with $M=dn$, $H_1=H$ and $H_2=\emptyset$, and note that $\tM=M-2|H|$, $\tM_2=M_2+O(d|H|)$, and thus, 
\[
\frac{\tM_2^2}{\tM^3}=\frac{(d-1)^2}{dn}\left(1+O\left(\frac{|H|}{dn}\right)\right),\quad \frac{\tM_2}{\tM^2}=\frac{d-1}{dn}\left(1+O\left(\frac{|H|}{dn}\right)\right),\quad\frac{\tM_2}{\tM}=(d-1)\left(1-O\left(\frac{|H|}{dn}\right)\right)
\]
and
\[
\sum_{(x,y)\in \W}\frac{\td_x\td_y}{\tM^2}=O\left(\frac{d^2|H|}{M^2}\right)=O\left(\frac{|H|}{n^2}\right).
\]
\fix{By letting $\bar\phi_H=\bar \phi_{H,\emptyset}$ and using the approximation of $\tM_2/\tM$ above, it is straightforward to find that $\bar\phi_H(uv)=\phi_H(uv)+d-d^2+2|H|+O(d|H|/n)$.}
Hence, 
\begin{align*}
\frac{1}{dn-2|H|}\left(1-\frac{\bar\phi_H(uv)}{dn-2|H|}\right)&=\frac{1}{dn}\left(1+\frac{2|H|}{dn}+O\left(\frac{|H|^2}{d^2n^2}\right)\right)\left(1-\frac{\bar\phi_H(uv)}{dn}+O\left(\frac{\bar\phi_H(uv)}{dn}\frac{|H|}{dn}\right)\right)\\
&=\frac{1}{dn}\left(1-\frac{\bar\phi_H(uv)-2|H|}{dn}+O\left(\frac{|H|^2}{d^2n^2}+\frac{|H|}{n^2}\right)\right)\\
&=\frac{1}{dn}\left(1-\frac{\phi_H(uv) \fix{+d-d^2+O(d|H|/n)} }{dn}+O\left(\frac{|H|^2}{d^2n^2}+\frac{|H|}{n^2}\right)\right)
\end{align*}
as $\bar\phi_H(uv)=O(d^2)$. Moreover,
\[
1-\frac{\tM_2^2}{\tM^3}-\frac{\tM_2}{\tM^2}-\sum_{(x,y)\in \W}\frac{\td_x\td_y}{\tM^2} = \left(1-\frac{d-1}{n} \right)\left(1+O\left(\frac{|H|}{n^2}\right)\right).
\]
Theorem~\ref{thm:special} follows by noting further that $d^5n/M^3=O(d^2/n^2)$ \fix{and
\[
\left(1-\frac{\phi_H(uv) +d-d^2}{dn}\right)\left(1-\frac{d-1}{n} \right)=\left(1-\frac{\phi_H(uv)}{dn}\right)\left(1+O\left(\frac{d^2}{n^2}\right) \right).\qed
\]

}

\subsection{Expectation and variance of strictly balanced subgraphs}
Let $H$ be a fixed graph. Let
\[
t=|V(H)| \quad \mbox{and} \quad h=|H|.
\]
Let $Z_{H}$ denote the number of subgraphs of $\G(n,d)$ that are isomorphic to $H$; each such subgraph is called a copy of $H$. Let 
\[
\mu_H=\ex Z_H, \quad \mbox{and}\quad \sigma^2_H=\var Z_H.
\]
Define, for each $1\le j\le h-1$,
\[
\rho_H(j)=\sup\left\{\rho: |V(H')|\ge t\left(\frac{|H|'}{|H|}+\rho\right)\ \mbox{for all subgraphs $H'$ of $H$ where $|H'|=j$}\right\}.
\]
See Example~\ref{example_variance} below for values of $\rho_H$ where $H$ is a cycle or  a clique.

Given a real number $x$ and a nonnegative integer $k$, let $(x)_k=\prod_{j=0}^{k-1}(x-j)$. Let $\text{aut}(H)$ denote the size of the automorphism group of $H$.

\begin{thm}\label{thm:variance}
Let $H$ be a fixed graph with $t$ vertices and $h$ edges. 
Then, the expectation of $Z_H$ is
\[
\mu_H=\frac{(n)_t}{\text{aut}(H)}\frac{\varphi(d,H)}{(dn)^{h}}\left(1-\frac{A_{H}d+B_{H}}{dn}+O(d^2/n^2)\right), 
\]
where 
\[
\varphi(d,H)=\prod_{v\in V(H)} (d)_{d_v^H},
\]
and
$A_H$ and $B_H$ are constants depending only on $H$.
Suppose further that $H$ is strictly balanced and
\[
\frac{d^{h-1}}{n^{h-t+1}}=o(1),\quad \frac{d^{h+2}}{n^{h-t+2}}=o(1),\quad \frac{\mu_H^{1-j/h}}{n^{t\rho_H(j)}}=o(1), \ \mbox{for every integer $1\le j\le h-1$}. 
\]
Then,
\[
\sigma_H^2\sim \mu_H \sim \frac{(n)_t}{\text{aut}(H)}\frac{\varphi(d,H)}{(dn)^{h}}.
\]
\end{thm}
\begin{example}\label{example_variance}
{\em (Cycles)} Let $H=C_{\ell}$ where $\ell\ge 3$. Then, $\rho_{C_{\ell}}(j)=1/\ell$ for every $1\le j\le \ell-1$. The assumptions of Theorem~\ref{thm:variance} are satisfied when $d^{\ell-1}=o(n)$ and $d^{\ell+2}=o(n^2)$. In particular, for $H=K_3$, the assumptions of Theorem~\ref{thm:variance} are satisfied when $d=o(n^{2/5})$.

{\em (Cliques)} Let $H=K_t$ where $t\ge 3$. Then, $\rho_{K_t}(j) = \frac{s(t-1)-2j}{t(t-1)}$ for every $\binom{s-1}{2}< j\le \binom{s}{2}$. The assumptions of Theorem~\ref{thm:variance} are equivalent to 
\[
d\ll n^{1-2/(t+1)},\quad d\ll n^{1-2t/(t^2-t+4)},\quad d\ll n^{(t+s-3)/(t+s-1)},
\]
which reduce to $d\ll n^{1-2t/(t^2-t+4)}$.

{\em (Trees)} Theorem~\ref{thm:variance} does not cover trees, as $d^{h-1}=o(n^{h-t+1})$ cannot be satisfied. Indeed the variance of $Z_{H}$ cannot be of order $\mu_H$ when $H$ is a tree. Consider $H=P_3$ to be a path of length three. Then $Z_H=d(d-1)^2n-6Z_{C_3}$, and therefore $\sigma_{P_3}^2=36\sigma_{C_3}^2\sim 36\mu_{C_3}$, which is much smaller than $\mu_{P_3}$.

\end{example}

\subsection{Distribution of the number of triangles}

\begin{thm}\label{thm:triangle}
Suppose $d=O(\sqrt{n})$ and $d=\omega(1)$. Let $Z_{C_3}$ be the number of triangles in $\G(n,d)$. Then,
\[
\frac{Z_{C_3}-\mu_{C_3}}{\sqrt{\mu_{C_3}}} \xrightarrow{d} \N(0,1),\quad \mbox{as $n\to \infty$}.
\]
\end{thm}

We will prove Theorem~\ref{thm:conditional} in Section~\ref{sec:conditional}. The variance of $Z_H$ for strictly balanced $H$ will be estimated in Section~\ref{sec:variance}. Finally, the limiting normal distribution of the number of triangles will be proved in Section~\ref{sec:distribution}. The proof for the case where $d=o(\sqrt{n})$ is given in Section~\ref{sec:small} and the case where $d=\Theta(\sqrt{n})$ is treated in Section~\ref{sec:large}. The proof for Corollary~\ref{cor:conditional} is presented in Section~\ref{sec:corollary}.

\section{Proof of Theorem~\ref{thm:conditional}}
\label{sec:conditional}

The proof is obtained by expressing $\pr(uv\in \G(n,\bfd)\mid H_1^+,H_2^-)$ by a function involving such conditional probabilities for other pairs $u'v'$. Applying Theorem~\ref{thm0:conditional} for pairs $u'v'$ yields a new estimate for  $\pr(uv\in \G(n,\bfd)\mid H_1^+,H_2^-)$ with improved error.
By repeatedly applying this argument it is possible to further improve the relative error and to relax conditions on $H_1$ and $H_2$. However, the expression of the conditional probability would become  more complicated after each iteration. We did not attempt this in the paper.

Let $\G$ be the class of graphs $G$ in $\G(n,\bfd)$ such that $H_1\subseteq G$ and $H_2\cap G=\emptyset$. Let $\G^{+}\subseteq \G$ be the set of graphs containing edge $uv$ and let $\G^{-}=\G\setminus\G^{+}$. We will estimate $|\G^+|/|\G^-|$ by defining switchings that relate graphs in $\G^+$ to graphs in $\G^-$ and counting the number of switchings that can be applied to a graph $G\in \G^+$, and the number of switchings that can produce a graph $G'\in \G^-$.

\begin{figure}
  \[
  \includegraphics[scale=1]{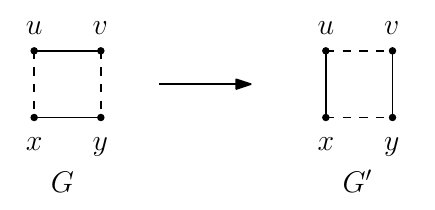}  
  \]
  \caption{Forward switching}
  \label{f:switch}
\end{figure}

Define the switching as follows.  
Given $G\in \G^+$, a forward switching specifies an ordered pair of vertices $(x,y)$ such that 
\begin{enumerate}
\item[(a)] $xy$ in $G\setminus H_1$ and $u$, $v$, $x$ and $y$ are all vertex distinct;
\item[(b)] none of $ux$ and $vy$ is  in $G\cup H_2$.
\end{enumerate} 
Then the forward switching replaces edges $uv$ and $xy$ by $ux$ and $vy$. The resulting graph $G'$ is obviously in $\G^-$. The inverse operation which converts $G'$ to $G$ is called a backward switching. See Figure~\ref{f:switch} for an illustration.

 Let $f(G)$ be the number of forward switchings that can be applied to $G$. The number of ways to choose $(x,y)$ such that $xy\in G\setminus H_1$ is $\tM$. Among these, the number of choices where $x\in\{u,v\}$ is $\td_u+\td_v$, and there are the same number of choices where $y\in\{u,v\}$. The number of choices where $\{x,y\}=\{u,v\}$ is 2. Thus, by inclusion-exclusion, the number of choices for $(x,y)$ satisfying (a) is $\tM-2(\td_u+\td_v)+2$.
To bound $f(G)$ we will subtract the number of choices for $(x,y)$ satisfying (a) but not (b).


Let $X_u$ be the number of choices of $(x,y)$ satisfying (a) where $ux$ is an edge in $G\cup H_2$,  and let $X_v$ be the number of choices of $(x,y)$ satisfying (a) where $vy$ is an edge in $G\cup H_2$. Let $X_{uv}$ be the number of choices  of $(x,y)$ satisfying (a) where both $ux$ and $vy$ are edges in $G\cup H_2$ and $xy\in G\setminus H_1$. Then,

\[
f(G)=\tM-2(\td_u+\td_v)+2-X_u-X_v+X_{uv}.
\]

\begin{claim} \label{claim:f} $\ex X_{uv}=O(\Delta^4/M)$ and
\[
\ex X_u=\sum_{x\in {\cal N}_{H_1\cup H_2}(u)} \td_x +\frac{(\td_u-1)\tM_2}{\tM}   +O(\Delta^4/M), \quad \ex X_v=\sum_{y\in {\cal N}_{H_1\cup H_2}(v)} \td_y +\frac{(\td_v-1)\tM_2}{\tM}   +O(\Delta^4/M).
\]
\end{claim}
By the claim,
\begin{align*}
\ex\, f(G)&=\tM-2(\td_u+\td_v)+2-\sum_{x\in {\cal N}_{H_1\cup H_2}(u)} \td_x - \sum_{y\in {\cal N}_{H_1\cup H_2}(v)} \td_y-\frac{(\td_u+\td_v-2)\tM_2}{\tM}  +O\left(\Delta^4/M\right)\\
&\hspace{-1cm}=\left(\tM-2(\td_u+\td_v)+2-\sum_{x\in {\cal N}_{H_1\cup H_2}(u)} \td_x - \sum_{y\in {\cal N}_{H_1\cup H_2}(v)} \td_y-\frac{(\td_u+\td_v-2)\tM_2}{\tM}\right)\left(1+O\left(\frac{\Delta^4}{M^2}\right)\right).
\end{align*}

Next, we estimate $b(G')$, the number of backward switchings that can be applied to $G'$. To perform a backward switching, we choose $x$ and $y$ such that $xu\in G'\setminus H_1$ and $yv\in G'\setminus H_1$ and $xy\notin G'\cup H_2$. Then replace $ux$, $yv$ by $uv$ and $xy$. The number of ways to choose $(x,y)$ such that $xu\in G'\setminus H_1$ and $yv\in G'\setminus H_1$ is $\td_u\td_v$. Thus, $b(G')=\td_u\td_v-Y_1-Y_2$ where $Y_1$ is the number of pairs $(x,y)$ such that $xu\in G'\setminus H_1$, $yv\in G'\setminus H_1$ and $xy\in G'\cup H_2$, and $Y_2$ is the number of vertices $x$ such that $xu\in G'\setminus H_1$ and $xv\in G'\setminus H_1$.  To obtain upper and lower bounds for $\ex b(G')$ it is sufficient to estimate $\ex Y_1$ and $\ex Y_2$. 

\fix{
\begin{claim}\label{claim:b}
\begin{eqnarray*}
\ex Y_1&=& \frac{\td_u\td_v \tM_2^2}{\tM^3} +\sum_{(x,y)\in \W} \frac{\td_x\td_y\td_u\td_v}{\tM^2}+O\left(\td_u\td_v\frac{\Delta^5n}{M^3}\right)\\
\ex Y_2&=&\frac{\td_u\td_v}{\tM^2}\tM_2+O\left(\td_u\td_v\frac{\Delta^3}{M^2}\right).
\end{eqnarray*}
\end{claim}
}
By the claim and noting that $\Delta^3/M^2=O(\Delta^5n/M^3)$ we have
\begin{eqnarray*}
\ex\, b(G')&=&\td_u\td_v- \frac{\td_u\td_v \tM_2^2}{\tM^3} -\frac{\td_u\td_v\tM_2}{\tM^2}-\sum_{(x,y)\in \W} \frac{\td_x\td_y\td_u\td_v}{\tM^2}+O\left(\td_u\td_v\frac{\Delta^5n}{M^3}\right)\\
&=&\td_u\td_v\left(1-\frac{\tM_2^2}{\tM^3}-\frac{\tM_2}{\tM^2}-\sum_{(x,y)\in \W} \frac{\td_x\td_y}{\tM^2}\right)\left(1+O\left(\frac{\Delta^5n}{M^3}\right)\right).
\end{eqnarray*}
\fix{
Note that 
$
\sum_{G\in \G^{+}} f(G) = \sum_{G'\in \G^{-}} b(G')
$
and thus
$
|\G^{+}| \cdot \ex f(G) = |\G^{-}| \cdot \ex b(G').
$
}
Therefore,
 noting that $\Delta^4/M^2=O(\frac{\Delta^5n}{M^3})$ we have
\begin{eqnarray}
&&\pr\Big(uv\in\G(n,d)\mid H_1^+, H_2^-\Big) =\frac{|\G^{+}|}{|\G^{+}|+|\G^{-}|} = \frac{\ex\, b}{\ex\, f +\ex\, b} \nonumber\\
&&\hspace{0.5cm}= \frac{\td_u\td_v\left(1-\frac{\tM_2^2}{\tM^3}-\frac{\tM_2}{\tM^2}-\sum_{(x,y)\in \W} \frac{\td_x\td_y}{\tM^2}\right) \left(1+O\left(\frac{\Delta^5n}{M^3}\right)\right)}{\tM-2(\td_u+\td_v)+2-\sum_{x\in {\cal N}_{H_1\cup H_2}(u)} \td_x-\sum_{y\in {\cal N}_{H_1\cup H_2}(v)} \td_y-\frac{(\td_u+\td_v-2)\tM_2}{\tM} +\td_u\td_v }\nonumber\\
&&\hspace {0.5cm}=\frac{\td_u\td_v}{\tM}\left(1-\frac{\tM_2^2}{\tM^3}-\frac{\tM_2}{\tM^2}-\sum_{(x,y)\in \W} \frac{\td_x\td_y}{\tM^2}\right)\left(1-\frac{\bar\phi_{H_1,H_2}(uv)}{dn-2|H_1|}\right) \left(1+O\left(\frac{\Delta^5n}{M^3}\right)\right).\nonumber\qed
\end{eqnarray}
\remove{
The theorem follows by noting that $\tM=M-2|H_1|$, $\tM_2=M_2-(2d-1)|H_1|-O(\Delta_{H_1}|H_1|)$, and thus,
\begin{eqnarray}
&&\frac{\tM_2^2}{\tM^3}=\frac{(d-1)^2}{dn}\left(1+\frac{(2d-4)|H_1|}{d(d-1)n}\right)\left(1+O\left(\frac{|H_1|^2}{d^2n^2}+\frac{\Delta_{H_1}|H_1|}{d^2n}\right)\right)\label{err1}\\
&&\frac{\tM_2}{\tM^2}=\frac{d-1}{dn}\left(1+\frac{(2d-3)|H_1|}{d(d-1)n}\right)\left(1+O\left(\frac{|H_1|^2}{d^2n^2}+\frac{\Delta_{H_1}|H_1|}{d^2n}\right)\right)\label{err2}\\\bar\phi_{H_1,H_2}
&&\frac{\tM_2}{\tM}=(d-1)\left(1-\frac{|H_1|}{d(d-1)n}\right)\left(1+O\left(\frac{|H_1|^2}{d^2n^2}+\frac{\Delta_{H_1}|H_1|}{d^2n}\right)\right)\label{err3}\\
&&1-\frac{-2(\td_u+\td_v)+2-\sum_{x\in {\cal N}_{H_1\cup H_2}(u)} \td_x-\sum_{y\in {\cal N}_{H_1\cup H_2}(v)} \td_y-\frac{(\td_u+\td_v-2)\tM_2}{\tM} +\td_u\td_v }{\tM}\nonumber\\
&&\hspace{0.7cm}=\left(1-\frac{\bar\phi_{H_1,H_2}(uv)}{dn-2|H_1|}\right)(1+O(\xi))
\end{eqnarray}
and noting that the additional errors caused by replacing $\tM_2^2/\tM^3$, $\tM_2/\tM^2$ and $\tM_2/\tM$ terms in the expression by the leading terms in~\eqn{err1}--\eqn{err3} respectively are all bounded by $O(|H_1|\Delta_{H_1}/dn^2)$, which, together with a few relative error of order at most $|H_1|/dn^2$, are all absorbed by $\xi$. \qed
}
Now we prove Claims~\ref{claim:f} and~\ref{claim:b}.

{\em Proof of Claim~\ref{claim:f}.\ }
First we estimate $\ex X_u$. Let $X_{u,1}$ be the number of $(x,y)$ such that $ux\in H_1\cup H_2$, $xy\in G\setminus H_1$ and $y\notin\{u,v\}$, $X_{u,2}$ be the number of $(x,y)$ such that $ux\in G\setminus H_1$ and $xy\in G\setminus H_1$ and $y\neq u$, and $X_{u,3}$ be the number of $(x,y)$ such that $ux\in G\setminus H_1$ and $xy\in G\setminus H_1$ and $y=v$. Then $X_u=X_{u,1}+X_{u,2}-X_{u,3}$. Obviously, 
\[
X_{u,1}=\sum_{x\in {\cal N}_{H_1\cup H_2}(u)} \td_x - X',
\]
where $X'$ is the number of $(x,y)$ such that $ux\in H_1\cup H_2$, $xy\in G\setminus H_1$ and $y=v$. Then, by Theorem~\ref{thm0:conditional},
\[
\ex X'= \sum_{\substack{x\in{\cal N}_{H_1\cup H_2}(u) \\ xv\notin H_1}} \pr(xv\mid (H_1\cup\{uv\})^+,H_2^-) = O(\Delta_{H_1\cup  H_2} \Delta^2/M)=O(\Delta^3/M).
\]
Thus, $\ex\, X_{u,1}=\sum_{x\in {\cal N}_{H_1\cup H_2}(u)} \td_x +O(\Delta^3/M)$.
On the other hand,
\begin{eqnarray*}
\ex X_{u,2}&=& \sum_{x\in [n]\setminus (\{u,v\} \cup {\cal N}_{H_1\cup H_2}(u))} \pr(ux \mid (H_1\cup\{uv\})^+,H_2^-) \cdot ({\td}_x-1).
\end{eqnarray*}
 Again by Theorem~\ref{thm0:conditional},
\begin{eqnarray*}
\ex X_{u,2}&=& \sum_{x\in [n]\setminus (\{u,v\} \cup {\cal N}_{H_1\cup H_2}(u))} \frac{(\td_u-1)\td_x}{\tM} (1+O(\Delta^2/M)) \cdot ({\td}_x-1)\\
&=&(\td_u-1) \sum_{x\in [n]\setminus (\{u,v\} \cup {\cal N}_{H_1\cup H_2}(u))} \frac{(\td_x)_2}{\tM}  +O(\Delta^4/M)\\
&=&\frac{\td_u-1}{\tM} \left(\tM_2- \sum_{z\in \{u,v\} \cup {\cal N}_{H_1\cup H_2}(u)} (\td_z)_2\right)  +O(\Delta^4/M)=\frac{(\td_u-1)\tM_2}{\tM}+O\left(\Delta^4/M\right).
\end{eqnarray*}
Finally, by Theorem~\ref{thm0:conditional}
\begin{eqnarray*}
\ex X_{u,3}&=& \sum_{x\in [n]\setminus (\{u,v\} \cup {\cal N}_{H_1\cup H_2}(u))} \pr(ux,vx \mid (H_1\cup\{uv\})^+,H_2^-)\\
&\le & (1+O(\Delta^2/M)) \sum_{x\in [n]}\frac{d_ud_v(d_x)_2}{\tM^2} =O(\Delta^3/M).
\end{eqnarray*}
Combining all together, we have
\[
\ex X_u=\sum_{x\in {\cal N}_{H_1\cup H_2}(u)} \td_x +\frac{(\td_u-1)\tM_2}{\tM}   +O\left(\frac{\Delta^4}{M}\right),
\]
as $\Delta_{H_1}\le \Delta$ and $\Delta_{H_2}=O(\Delta)$ by assumption.
By symmetry,
\[
\ex X_v=\sum_{y\in {\cal N}_{H_1\cup H_2}(v)} \td_y+\frac{(\td_v-1)\tM_2}{\tM}  +O\left(\frac{\Delta^4}{M}\right).
\]
Finally, we bound $\ex X_{uv}$ by $O(\Delta^4/M)$. First we expose the neighbours of $u$ and $v$ in $G$. Next, for each $x\in {\cal N}_{G\cup H_2}(u)$ and $y\in {\cal N}_{G\cup H_2}(v)$, if $x\neq y$ and $xy\notin H_1$, the probability that $xy\in \G(n,\bfd)$, conditional on $H_1^+$, $H_2^-$, and the presence of the exposed edges incident with $u$ and $v$,  is at most $(1+O(\Delta^2/M))d^2/(\tM-4\Delta)=O(\Delta^2/M)$. There are at most $O(\Delta^2)$ such pairs of $(x,y)$, as $\Delta_{H_2}=O(\Delta)$. Thus,
\begin{eqnarray*}
\ex X_{uv}= O(\Delta^2) \cdot O(\Delta^2/M)=O(\Delta^4/M).  \qed
\end{eqnarray*}

{\em Proof of Claim~\ref{claim:b}.\ }
Let $Y_{1,1}$ be the number of pairs $(x,y)$ such that $xu\in G'\setminus H_1$ and $yv\in G'\setminus H_1$ and $xy\in H_1\cup H_2$ and $Y_{1,2}$ be the number of pairs $(x,y)$ such that $xu\in G'\setminus H_1$ and $yv\in G'\setminus H_1$ and $xy\in G'\setminus H_1$. Then $Y_1=Y_{1,1}+Y_{1,2}$. We first estimate $\ex Y_{1,1}$. Recall that $\W$ is the set of $(x,y)$ such that $xy\in H_1\cup H_2$, $xu,yv\notin H_1\cup H_2$, and $\{x,y\}\cap\{u,v\}=\emptyset$. 
 
%
Now, 
 \[
\ex Y_{1,1}=\sum_{(x,y)\in \W} \pr(ux,vy\in\G(n,d)\mid H_1^+, (H_2\cup\{uv\})^-).
\]
By Theorem~\ref{thm0:conditional},
\begin{eqnarray}
\ex Y_{1,1}&=&(1+O(\Delta^2/M)){\sum_{(x,y)\in \W}}\frac{\td_u\td_x}{\tM} \frac{\td_y\td_v}{\tM}={\sum_{(x,y)\in \W}}\frac{\td_u\td_x}{\tM} \frac{\td_y\td_v}{\tM} +O\left(\td_u\td_v\frac{\Delta^5n}{M^3}\right),
\end{eqnarray}
since $|\W|=O(\Delta n)$.
Similarly,
\[
\ex Y_{1,2}={\sum} \frac{\td_u(\td_x)_2(\td_y)_2\td_v}{\tM^3}\left(1+O\left(\frac{\Delta^2}{M}\right)\right)={\sum}\frac{\td_u(\td_x)_2(\td_y)_2\td_v}{\tM^3}+O\left(\td_u\td_v\frac{\Delta^5n}{M^3}\right),
\]
where the summation is over all pairs $(x,y)$ such that $x\in[n]\setminus(\{u,v\}\cup {\cal N}_{H_1\cup H_2}(u))$ and $y\in [n]\setminus(\{u,v,x\}\cup {\cal N}_{H_1\cup H_2}(v)\cup {\cal N}_{H_1\cup H_2}(x))$. Thus,
\[
\ex Y_{1,2}=O\left(\td_u\td_v\frac{\Delta^5n}{M^3}\right)+ \frac{\td_u\td_v}{\tM^3}\left(\tM_2^2 -a_1-a_2 \right),
\]
where 
\begin{eqnarray*}
a_1&=&\tM_2\left(\sum_{x\in \{u,v\}\cup {\cal N}_{H_1\cup H_2}(u)} (\td_x)_2 + \sum_{j\in \{u,v\}\cup {\cal N}_{H_1\cup H_2}(v)} (\td_y)_2\right)\\
&&\hspace{1cm}-\sum_{x\in \{u,v\}\cup {\cal N}_{H_1\cup H_2}(u)} (\td_x)_2 \sum_{j\in \{u,v\}\cup {\cal N}_{H_1\cup H_2}(v)} (\td_y)_2= O(\Delta^5 n),
\end{eqnarray*}
and
$
a_2= \sum (\td_x)_2(\td_y)_2,
$
where the sum is over all pairs $(x,y)$ such that $y\in \{x\}\cup {\cal N}_{H_1\cup H_2}(x)$ in addition to $x\in[n]\setminus(\{u,v\}\cup {\cal N}_{H_1\cup H_2}(u))$ and $y\in[n]\setminus(\{u,v\}\cup {\cal N}_{H_1\cup H_2}(v))$. Trivially, $a_2=O(\Delta^4|H_1\cup H_2|)=O(\Delta^5n)$.
Thus,
\[
\ex Y_{1,2}= \frac{\td_u\td_v\tM_2^2}{\tM^3}  +O\left(\td_u\td_v\frac{\Delta^5n}{M^3}\right).
\]
It follows now that
\[
\ex Y_1= \frac{\td_u\td_v \tM_2^2}{\tM^3} +{\sum_{(x,y)\in \W}}\frac{\td_u\td_v\td_x\td_y}{\tM^2} +O\left(\td_u\td_v\frac{\Delta^5n}{M^3}\right).
\]
With a similar but simpler argument, we have 
\begin{eqnarray*}
\ex Y_2&=&(1+O(\Delta^2/M))\sum_{x\in [n]\setminus (\{u,v\}\cup {\cal N}_{H_1\cup H_2}(u)\cup {\cal N}_{H_1\cup H_2}(v))} \frac{\td_u\td_v(\td_x)_2}{\tM^2}\\
&=&(1+O(\Delta^2/M))\frac{\td_u\td_v}{\tM^2}\Big(\tM_2-O(\Delta^3)\Big)= \frac{\td_u\td_v\tM_2}{\tM^2} +O\left(\td_u\td_v\frac{\Delta^3}{M^2}\right). \qed
\end{eqnarray*}
\remove{
Given a degree sequence $\bfd$, let $\Delta(\bfd)$ denote the maximum component in $\bfd$ and $\delta(\bfd)$ denote the minimum component in $\bfd$, and let $\text{rng}(\bfd)=\Delta(\bfd)-\delta(\bfd)$.

\begin{thm}\label{thm:nonconstant} Let $\bfd$ be a degree sequence with $\Delta(\bfd)=O(\delta(\bfd))$. Let $\Delta=\Delta(\bfd)$.
Let $H_1$ and $H_2$ be disjoint graphs with degree sequences $\bfd^{H_1}\preceq \bfd$ and $\Delta_{H_2}=O(\Delta)$. Let $\tbfd=\bfd-\bfd^{H_1}$. Suppose $uv\notin H_1\cup H_2$, and let 
\[
\xi=\frac{\Delta^4}{\td_u\td_vn^2}+|H_1\cup H_2|\left(\frac{(\Delta_{H_1\cup H_2} +\text{rng}(\bfd)) \Delta}{\td_u\td_vn^2}+\frac{\Delta^3}{\td_u\td_vn^3}+\frac{\Delta|H_1|}{\td_u\td_vn^3}\right).
\]
Then $\pr\Big(uv\in\G(n,d)\mid H_1^+, H_2^-\Big)$ is given by~\eqn{probability} with the term $d^2|H_1\cup H_2|/\td_u\td_vn^2$ replaced by $\Delta^2|H_1\cup H_2|/\td_u\td_vn^2$.
\end{thm}

\proof The proof is exactly the same as for Theorem~\ref{thm:conditional}, except for the following two changes.
\begin{itemize}
\item The error $O(d/n)$ becomes $O(\Delta/n)$ when Theorem~\ref{thm0:conditional} is applied. All error terms in the proof of Theorem~\ref{thm:conditional} remain the same after replacing $d$ by $\Delta$.

\item The term $d^2|H_1\cup H_2|/\td_u\td_vn^2$ will be replaced by $\Delta^2|H_1\cup H_2|/\td_u\td_vn^2$ in Claim~\ref{claim:b}. Now all vertices in $\G(n,\bfd)$ has degree $\Delta+O(\text{rng}(\bfd))$. Thus,~\eqn{Y11} changes to
\begin{eqnarray}
\ex Y_{1,1}&=&(1+O(\Delta/n)){\sum_{\substack{(x,y)\notin \W\\ xy\in  H_1\cup H_2}}}\frac{\td_u\td_x}{\tM} \frac{\td_y\td_v}{\tM}= (1+O(\Delta/n)){\sum_{\substack{(x,y)\notin \W\\ xy\in  H_1\cup H_2}}}\frac{(\Delta-O(\text{rng}(\bfd)+\Delta_{H_1}))^4}{(M-2|H_1|)^2}\nonumber\\
  &=&\left((1+O\left(\frac{\Delta}{n}+\frac{|H_1|}{\Delta n}\right)\right)\left(|H_1\cup H_2|-O(\Delta^2)\right)\frac{\Delta^4+O(\Delta^3(\text{rng}(\bfd)+\Delta_{H_1}))}{d^2n^2}\nonumber\\
  &=&\frac{\Delta^2|H_1\cup H_2|}{n^2}+O\left(\frac{\Delta(\Delta_{H_1}+\text{rng}(\bfd))|H_1\cup H_2|}{n^2}+\frac{\Delta^4}{n^2}+|H_1\cup H_2|\left(\frac{\Delta^3}{n^3}+\frac{\Delta|H_1|}{n^3}\right)\right).\nonumber
\end{eqnarray}
\end{itemize}
These changes explain the replacement of $d^2|H_1\cup H_2|/\td_u\td_vn^2$ by $\Delta^2|H_1\cup H_2|/\td_u\td_vn^2$ in the expression for  $\pr\Big(uv\in\G(n,d)\mid H_1^+, H_2^-\Big)$. They also explain the modifications, in particular the additional term $\text{rng}(\bfd)$, in $\xi$.\qed

}


\section{Variance of the number of strictly balanced subgraphs}
\label{sec:variance}

Let $H$ be a fixed graph with $t$ vertices and $h$ edges.
Recall that $Z_{H}$ denotes the number of subgraphs of $\G(n,d)$ that are isomorphic to $H$, and
$\mu_H=\ex Z_H$, $\sigma^2_H=\var Z_H$.
The goal of this section is to estimate $\sigma^2_H$ for strictly balanced $H$. We start by estimating $\mu_H$.

\subsection{Expectation}

Given a graph $H$, let 
\[
\varphi_H = \prod_{v\in V(H)} (d)_{d^H_v}.
\]
Let $J$ be a copy of $H$ in $K_n$, and let 
\[
\lambda_H = \pr(J\subseteq \G(n,d)).
\]  
By Theorem~\ref{thm:special}, 
\begin{equation}
\lambda_H=\frac{\varphi(d,H)}{(dn)^{h}}\left(1-\frac{A_{H}d+B_{H}}{dn}+O(d^2/n^2)\right), \label{lambda_H}
\end{equation}
where $A_H$ and $B_H$ are constants that depend only on $H$. The values of $A_{H}$ and $B_H$ are given for $H$ being a cycle in  Example~\ref{triangle_expectation} below.
By linearity of expectation,
\begin{equation}
\mu_H=\frac{(n)_t}{\text{aut}(H)} \lambda_H, \label{mu_H}
\end{equation}
where $\text{aut}(H)$ denotes the size of the automorphism group of $H$.
Hence,
\[
\mu_H=\Theta(n^t (d/n)^{h}) = \Theta(d^{h} n^{t-h}).
\]

\begin{example}\label{triangle_expectation}
Let $H=C_{\ell}$, a cycle of length $\ell\ge 3$. By Theorem~\ref{thm:special}, 
\[
\lambda_{C_{\ell}}=\left(\frac{d-1}{n}\right)^{\ell} \left(1-\sum_{i=1}^{\ell}\frac{\phi_i}{dn}+O(d^2/n^2)\right), 
\] 
where
\[
\phi_1=-d,\quad \phi_2=-2d, \quad, \phi_i=-2d-2(i-3)-1\ \mbox{for all $3\le i\le \ell-1$},\quad \phi_{\ell}=-3d+3-2(\ell-3).
\]
In particular, 
\begin{equation}
\lambda_{C_{3}}=\left(\frac{d-1}{n}\right)^{3}\left(1+\frac{6d-3}{dn}+O(d^2/n^2)\right).\label{lambda_C3}
\end{equation}
\end{example}
\subsection{Variance of $Z_{C_3}$}
\label{sec:variance-triangle}
In this section, we calculate  $\sigma^2_{C_3}$ in detail. In the next section, we will briefly sketch how to extend the argument for $C_3$ to other strictly balanced graphs. Let $N= \binom{n}{3}$, and let $H_1,\ldots, H_N$ be an enumeration of all copies of triangles in $K_n$. Let $X_i=\ind{H_i\subseteq \G(n,d)}$. 
Then $Z_{C_3}=\sum_{i=1}^N X_i$. Moreover, $\ex X_i=\lambda$ for every $i$ where $\lambda=\lambda_{C_3}$ has been estimated in~\eqn{lambda_C3}. Hence,
\begin{equation}
\sigma^2_{C_3}=\ex (Z_{C_3}-\mu_{C_3})^2 = \ex\left(\sum_{i=1}^{N}(X_i-\lambda)\right)^2 =  \sum_{(i,j)\in [N]^2} (\ex X_iX_j -\lambda^2).\label{sigma_C3}
\end{equation}
We split the above sum into four parts according to how $H_i$ and $H_j$ intersect.

\begin{itemize}
\item[(a)]  $H_i$ and $H_j$ are vertex disjoint.
For such pairs $(i,j)$, by Theorem~\ref{thm:special} and~\eqn{phi1},
\[
\ex X_iX_j=\left(\frac{d-1}{n}\right)^{6}\left(1+\frac{12d+12}{dn}+O(d^2/n^2)\right).
\]
Thus,
\[
\ex X_iX_j-\lambda^2=\left(\frac{d-1}{n}\right)^{6}\frac{18}{dn}+O(d^8/n^8).
\]
The number of vertex disjoint pairs $(H_i, H_j)$ is
$\binom{n}{3}\binom{n-3}{3}=(n)_6/36$.
Hence, the contribution to~\eqn{sigma_C3} from part (a) is
\begin{equation}
\frac{(n)_6}{36}\left(\left(\frac{d-1}{n}\right)^{6}\frac{18}{dn}+O(d^8/n^8)\right)=\left(\frac{d-1}{n}\right)^{6}\frac{(n)_6}{2dn}+O(d^8/n^2). \label{sum1}
\end{equation}

\item[(b)] $H_i$ and $H_j$ share exactly one vertex. For such pairs,
\begin{eqnarray*}
\ex X_iX_j&=&\frac{(d(d-1))^4(d)_4}{(dn)^6}\left(1+\frac{16d-4}{dn}+O\left(d^2/n^2\right)\right)\\
&=&\frac{(d-1)^5(d-2)(d-3)}{dn^6}\left(1+\frac{16d-4}{dn}+O\left(d^2/n^2\right)\right).
\end{eqnarray*}
Thus,
\[
\ex X_iX_j-\lambda^2=-\frac{2(2d-3)(d-1)^5}{dn^6}+O(d^6/n^7+d^8/n^8).
\]
The number of pairs $(H_i,H_j)$ sharing exactly one vertex is
$\binom{n}{3} 3 \binom{n-3}{2}=(n)_5/4$.
Thus, the contribution to~\eqn{sigma_C3} from part (b) is
\begin{equation}
 \frac{(n)_5}{4}\left(-\frac{2(2d-3)(d-1)^5}{dn^6}+O\left(\frac{d^6}{n^7}+\frac{d^8}{n^8}\right)\right)=-\frac{(n)_5(2d-3)(d-1)^5}{2dn^6}+O\left(\frac{d^6}{n^2}+\frac{d^8}{n^3}\right).\label{sum2}
\end{equation}

\item[(c)] $H_i$ and $H_j$ share exactly two vertices.
For such pairs,
\[
\ex X_iX_j=\frac{(d)_2^2(d)_3^2}{(dn)^5}\left(1+\frac{13d-11}{dn}\right)=\frac{(d-1)^4(d-2)^2}{dn^5}\left(1+\frac{13d-11}{dn}\right).
\]
Thus,
\[
\ex X_iX_j-\lambda^2=\frac{(d-1)^4(d-2)^2}{dn^5}+O(d^6/n^6).
\]
The number of such pairs $(H_i, H_j)$ is
$\binom{n}{3}3(n-3)=(n)_4/2$. Thus, the contribution to~\eqn{sigma_C3} from part (c) is
\begin{equation}
 \frac{(n)_4}{2}\left(\frac{(d-1)^4(d-2)^2}{dn^5}+O(d^6/n^6)\right)=\frac{(n)_4(d-1)^4(d-2)^2}{2dn^5}+O(d^6/n^2). \label{sum3}
\end{equation}

\item[(d)]
$H_i=H_j$. For such pairs, \[
\ex X_iX_j=\ex X_i=\lambda.
\]
Thus,
\[
\ex X_iX_j-\lambda^2=\lambda-\lambda^2=\left(\frac{d-1}{n}\right)^{3}\left(1+\frac{6d-3}{dn}+O(d^2/n^2)\right)+O(d^6/n^6).
\]
There are $(n)_3/6$ pairs $(H_i,H_j)$ where $H_i=H_j$.
Thus, the contribution to~\eqn{sigma_C3} from part (d) is
\begin{equation}
 \frac{(n)_3}{6}\left(\left(\frac{d-1}{n}\right)^{3}\left(1+\frac{6d-3}{dn}\right)+O(d^5/n^5)\right)= \frac{(n)_3}{6}\left(\frac{d-1}{n}\right)^{3}\left(1+\frac{6d-3}{dn}\right)+O(d^5/n^2).\label{sum4}
\end{equation}
\end{itemize}
By~\eqn{sum1}--\eqn{sum3}, the total contribution to~\eqn{sigma_C3} from parts (a)--(c) is
\[
-\frac{(d-1)^4(d-2)(n-1)_3}{2dn^4}+O(d^8/n^2).
\]
Combining with~\eqn{sum4} and by~\eqn{sigma_C3} we have
\[
\sigma^2_{C_3}=\frac{(n)_3}{6}\left(\frac{d-1}{n}\right)^{3}\left(1+\frac{6d-3}{dn}\right)-\frac{(d-1)^4(d-2)(n-1)_3}{2dn^4}+O(d^8/n^2)=(1+o(1))\frac{(d-1)^3}{6}+O(d^8/n^2).
\]
Since
$
d=o(n^{2/5})
$, we have
\[
d^8/n^2=o(d^3).
\]
It follows now that
\[
\sigma^2_{C_3}\sim \mu_{C_3}\sim \frac{(d-1)^3}{6}.
\]

\subsection{Variance of $Z_H$: proof of Theorem~\ref{thm:variance}}

Recall that
\[
\rho_H(j)=\sup\left\{\rho: |V(H')|\ge t\left(\frac{|H|'}{|H|}+\rho\right)\ \mbox{for all subgraph $H'\neq H$ of $H$ where $|H'|=j$}\right\},
\]
and $\eta_H(j)$ is positive for every $1\le j\le h-1$ as $H$ is strictly balanced.

Recall $\lambda_H$ and $\mu_H$ from~\eqn{lambda_H} and~\eqn{mu_H}. 
Let $N= \binom{(n)_t}{\text{aut}(H)}$, and let $H_1,\ldots, H_N$ be an enumeration of all copies of $H$ in $K_n$. Let $X_i=\ind{H_i\subseteq \G(n,d)}$. Then, $Z_{H}=\sum_{i=1}^N X_i$. Moreover, $\ex X_i=\lambda$ for every $i$ where $\lambda=\lambda_{H}$. Hence,
\begin{equation}
\sigma^2_{H}=\ex (Z_{H}-\mu_{H})^2 = \ex\left(\sum_{i=1}^{N}(X_i-\lambda)\right)^2 =  \sum_{(i,j)\in [N]^2} (\ex X_iX_j -\lambda^2).\label{sigma_H}
\end{equation}
Similarly to the case $H=C_3$, we split the above sum into four parts according to the intersection of $H_i$ and $H_j$. We will give only upper bounds for the contributions from the first three parts, and estimate asymptotically the contribution from the last part. Under the hypotheses of the theorem, we verify that the contribution from the last part dominates. The calculations are similar to the case where $H=C_3$ and thus we only briefly sketch the proof.
\begin{itemize}
\item[(a)] $H_i$ and $H_j$ are vertex disjoint. There are $O(n^{2t})$ such pairs. For each such pair,  
\[
\ex X_iX_j-\lambda^2=O(1/dn+d^2/n^2) \left(\frac{d}{n}\right)^{2h}.
\]
Hence the total contribution to~\eqn{sigma_H} from this part is $\mu_H\cdot O(d^{h-1}n^{t-h-1}+d^{h+2}n^{t-h-2})$.
\item[(b)] $H_i$ and $H_j$ are edge disjoint but share at least one vertex. There are $O(n^{2t-1})$ such pairs. For each such pair, by Theorem~\ref{thm:special} and~\eqn{phi2},
\begin{eqnarray*}
&& \ex X_iX_j=\lambda \cdot \frac{\prod_{v\in H}(d-O(1))_{d_v^H}}{(dn)^{h}} \left(1-\frac{O(d)}{dn}+O(d^2/n^2)\right) =\lambda^2(1+O(1/d+d^2/n^2)).
\end{eqnarray*}
Thus,
\begin{eqnarray*}
&& \ex X_iX_j-\lambda^2=\lambda^2\cdot O( 1/d+d^2/n^2) =\lambda \cdot O(d^{h-1}/n^h+d^{h+2}/n^{h+2}),
\end{eqnarray*}
and the total contribution to~\eqn{sigma_H} from this part is $\mu_H\cdot O(d^{h-1}n^{t-h-1}+d^{h+2}n^{t-h-3})$.

\item[(c)] $H_i$ and $H_j$ share $x\ge 2$ vertices and $1\le y\le h-1$ edges. There are $O(n^{2t-x})$ such pairs. 
For each such pair, we know $x\ge t(y/h+\rho_H(y))$ by definition of $\rho_H$. Thus,
\[
\ex X_iX_j-\lambda^2= O(\lambda^2+\ex X_iX_j)=  O(\lambda^2+(d/n)^{2h-y}).
\]
The total contribution to~\eqn{sigma_H} from part (c) is
\[
\sum_{x,y}\Big(O(\lambda^2 n^{2t-x}) +O(n^{2t-x} (d/n)^{2h-y}) \Big)= \mu_H \sum_{x,y} O(n^{t-x} (d/n)^{h-y}) = \mu_H \sum_{x,y}O(n^{-t\rho_H(y)} \mu_H^{1-y/h}).
\]

\item[(d)] $H_i=H_j$. The total contribution to~\eqn{sigma_H} from this part is
\[
\frac{(n)_t}{\text{aut}(H)} (\lambda-\lambda^2) \sim \mu_H\sim \frac{n^t}{\text{aut}(H)}\frac{\varphi(d,H)}{(dn)^h}.
\]
\end{itemize}

By the hypotheses of the theorem, the total contributions to~\eqn{sigma_H} from parts (a)--(c) are $o(\mu_H)$, and the assertion of the theorem follows.\qed

\section{Distribution of the number of triangles}
\label{sec:distribution}

We use the following result~\cite{gao2004asymptotic} by Z. Gao and Wormald to determine the limiting distribution of $Z_{C_3}$.
\begin{thm} \label{thm:normality}
Let $s_n>\mu_n^{-1}$ and $\sigma_n=\sqrt{\mu_n+\mu_n^2 s_n}$, where $0<\mu_n\to\infty$. Suppose that $\mu_n=o(\sigma_n^3)$, and a sequence $(X_n)$ of nonnegative random variables satisfies
\[
\ex (X_n)_k\sim \mu_n^k\exp\left(\frac{k^2s_n}{2}\right)
\]
uniformly for all integers $k$ in the range $c\mu_n/\sigma_n\le k\le c'\mu_n/\sigma_n$ for some constants $c'>c>0$. Then,
\[
\frac{X_n-\mu_n}{\sigma_n} \xrightarrow{d} \N(0,1),\quad \mbox{as $n\to\infty$.}
\]
\end{thm}

\subsection{Proof of Corollary~\ref{cor:conditional}}
\label{sec:corollary}
Let $e_1,\ldots, e_{|F|}$ be an enumeration of the edges in $F$. Let $F_0=\emptyset$, $H_0=H$, $F_j=F_{j-1}\cup\{e_j\}$ and $H_j=H\cup F_j$ for every $1\le j\le |F|$. Fix $1\le j\le |F|$. let $u$ and $v$ be the ends of $e_j$ and let $\tbfd=\bfd-\bfd^{H_{j-1}}$. By~\eqn{phi0} and~\eqn{phi2},
\begin{eqnarray*}
&&\pr\Big(e_j\in\G(n,d)\mid H_{j-1}^+\Big)=\frac{\td_u\td_v}{dn}\left(1-\frac{\phi_{H_{j-1}}(uv)}{dn}\right)\left(1+O\left(\frac{|H_{j-1}|}{n^2}+\frac{|H_{j-1}|^2}{d^2n^2}+\frac{d^2}{n^2}\right)\right)\\
&&=\left(\frac{\td_u\td_v}{dn}+O\left(\frac{d|H_{j-1}|}{n^3}+\frac{|H_{j-1}|^2}{dn^3}+\frac{d^3}{n^3}\right)\right)\left(1+\frac{d(1+d_u^{H_{j-1}}+d_v^{H_{j-1}})+O(|H_{j-1}|)}{dn}\right)\\
&&=\left(\frac{(d-d^{H_{j-1}}_u)(d-d_v^{H_{j-1}})}{dn}\right)\left(1+\frac{d(1+d_u^{H_{j-1}}+d_v^{H_{j-1}})}{dn}\right)+O\left(\frac{d^3}{n^3}+\frac{|H_{j-1}|}{n^2}\right),
\end{eqnarray*}
where the last step above holds since \fix{the errors $d|H_{j-1}/n^3|$ and $|H_{j-1}|^2/dn^3$ are absorbed by $|H_{j-1}|/n^2$ as} $|H_{j-1}|=O(dn)$.
It is easy to verify using elementary calculus that $$\left(\frac{(d-d^{H_{j-1}}_u)(d-d_v^{H_{j-1}})}{dn}\right)\left(1+\frac{d(1+d_u^{H_{j-1}}+d_v^{H_{j-1}})}{dn}\right)$$
is a decreasing function of $d^{H_{j-1}}_u$ and $d^{H_{j-1}}_v$, and thus the above product is at most
\[
\frac{(d-d^{F_{j-1}}_u)(d-d_v^{F_{j-1}})}{dn}\left(1+\frac{d(1+d_u^{F_{j-1}}+d_v^{F_{j-1}})}{dn}\right)=\frac{(d-d^{F_{j-1}}_u)(d-d_v^{F_{j-1}})}{dn}\left(1+O(n^{-1})\right).
\]
It follows that
\begin{eqnarray*}
\pr\Big(e_j\in\G(n,d)\mid H_{j-1}^+\Big)&\le& \frac{(d-d^{F_{j-1}}_u)(d-d_v^{F_{j-1}})}{dn}\left(1+O(n^{-1})\right)+O\left(\frac{d^3}{n^3}+\frac{|H_{j-1}|}{n^2}\right)\\
&=&\frac{(d-d^{F_{j-1}}_u)(d-d_v^{F_{j-1}})}{dn}\left(1+O\left(\frac{1}{n}+\frac{d^2}{n^2}+\frac{|H_{j-1}|}{dn}\right)\right).
\end{eqnarray*}
The assertion of the corollary follows by applying the above bound to
\[
\pr(F\subseteq \G(n,d)\mid H^+)=\prod_{j=1}^{|F|}\pr\Big(e_j\in\G(n,d)\mid H_{j-1}^+\Big)
\]
and then comparing the resulting expression with~\eqn{lambda_H}.\qed

\subsection{Distribution of $Z_{C_3}$ when $d=o(\sqrt{n})$}
\label{sec:small}

Let  $Z=Z_{C_3}$, $\lambda=\lambda_{C_3}$ and $\mu=\mu_{C_3}$. Then, $\mu=\Theta(d^3)$.
Assume $d=o(\sqrt{n})$ and $d\to\infty$. We aim to show that
\[
\ex (Z)_k \sim \mu^k,
\]
for all $k=O(d^{3/2})$. Then the distribution of $Z$ follows by Theorem~\ref{thm:normality} with $\mu_n=\mu$, $\sigma_n=\sqrt{\mu}$ and $s_n=0$. 

Let $N=\binom{n}{3}$ and let $H_1,\ldots, H_N$ be an enumeration of all copies of triangles in $K_n$. Let $\F$ be the set of  ordered $k$-tuples $(j_1, \ldots, j_k)\in [N]^k$ that are pairwise distinct. Then, for any positive integer $k$,
\begin{equation}
\ex (Z)_k = \sum_{(j_1,\ldots, j_k)\in\F} \pr\left(\cup_{i=1}^k H_{j_i}\subseteq \G(n,d)\right). \label{factorial}
\end{equation}

\begin{definition}\label{def:hole}
Given a $k$-tuple $\bfj$, we say a triple of vertices $\{x,y,z\}$ is a {\em hole} induced by ${j_1},\ldots, {j_{k}}$, if there exist ${j_{i_1}}$, ${j_{i_2}}$ and ${j_{i_3}}$, such that $H_{j_{i_1}}$ contains $xy$, $H_{j_{i_2}}$ contains $yz$ and $H_{j_{i_3}}$ contains $xz$, and none of $H_{j_{i_1}}$, $H_{j_{i_2}}$ and $H_{j_{i_3}}$ is $xyz$. 
\end{definition}
We partition the summation into the following two parts.
\begin{itemize}
\item[(a)] $\F_1\subseteq \F$:  the set of $(j_1, \ldots, j_k)$ such that $H_{j_1},\ldots, H_{j_k}$ are pairwise edge disjoint.
\item[(b)] $\F_2=\F\setminus\F_1$. 
\end{itemize}

We will prove the following.
\begin{lemma}\label{lem:F1}
\[
\sum_{(j_1,\ldots, j_k)\in\F_1} \pr\left(\cup_{i=1}^k H_{j_i}\subseteq \G(n,d)\right) \sim \mu^k.
\]
\end{lemma}

\begin{lemma}\label{lem:F2}
\[
\sum_{(j_1,\ldots, j_k)\in\F_2} \pr\left(\cup_{i=1}^k H_{j_i}\subseteq \G(n,d)\right) = o(\mu^k). 
\]
\end{lemma}


Now Theorem~\ref{thm:triangle} follows by the above two lemmas and Theorem~\ref{thm:normality}. It only remains to prove Lemmas~\ref{lem:F1} and~\ref{lem:F2}.

{\em Proof of Lemma~\ref{lem:F1}.\ } We say $H_{j_i}$ hits $v$ if $v$ is incident with an edge in $H_{j_i}$.
Given a $k$-tuple $H_{j_1},\ldots, H_{j_k}$, let $h(v)$ be the number of $H_j$ which hit $v$. Let $h_j$ be the number of vertices $v$ with $h(v)=j$. Then, given any $(j_1,\ldots, j_k)\in \F_1$,
the number of vertices in $H^{\oplus}=\cup_{i=1}^k H_i$ is equal to $vk-\sum_{j\ge 2}(j-1) h_j$ and the number of edges in $H^{\oplus}$ is $3k$.

Let $\B$ be the set of $(j_1,\ldots, j_k)\in [N]^k$ satisfying the following property.
\[
h_j>0 \ \mbox{for some $j\ge 4$; or}\ h_j \ge d^2/\sqrt{n}\ \mbox{for some $j\in \{2,3\}$}.
\]
We further partition $\F_1$ into two parts: $\F_1\setminus \B$ and $\F_1\cap \B$. 

\begin{claim}\label{claim:F1}
\[
|\F_1\setminus \B| \sim \left(\frac{(n)_3}{6}\right)^k,\ \mbox{and}\ |\F_1\cap \B|=o\left(\left(\frac{(n)_3}{6}\right)^k\right).
\]
\end{claim}

Fix an arbitrary $(j_1,\ldots, j_k)\in \F_1\setminus \B$. Let $e_1,\ldots, e_{3k}$ be an enumeration of the edges in $H^{\oplus}=\cup_{i=1}^k H_{j_i}$, where $e_{3\ell-2},e_{3\ell-1},e_{3\ell}$ is the set of edges in $H_{j_{\ell}}$. Let $F_i=\cup_{j<i} e_j$. By~\eqn{phi0} and~\eqn{phi2},
$\phi_{F_i}(e_i)=O(i+d)$. Thus, 
\[
\sum_{i=1}^k \frac{\phi_{F_i}(e_i)}{dn} = \frac{O(k^2+dk)}{dn}.
\] 
Since $\bfj\notin \B$, we have $d_u^{H^{\oplus}}=O(1)$ for every $u$. Hence, by Theorem~\ref{thm:special} (with $\td_u,\td_v=d-O(1)$ and $|H|=O(k)$)
\begin{equation}
\ex \left(\prod_{i=1}^k X_{H_{j_i}}\right) = \lambda^k\left(1+O\left(\frac{h_2+h_3}{d}\right)\right)\left(1+O\left(\frac{k^2+dk}{dn}+\frac{k^2+d^2k}{n^2}+\frac{k^3}{d^2n^2}\right)\right) \sim \lambda^k, \label{kthMoments}
\end{equation}
where the error $O((h_2+h_3)/d)$ accounts for the case that some $H_{j_{\ell}}$ lands on a vertex which has been occupied by some $H_{j_i}$ where $i<\ell$, and the last step of asymptotic relation above holds  
as $k=O(d^{3/2})$ and $d=o(\sqrt{n})$.  By Claim~\ref{claim:F1}, 
\begin{equation}
\sum_{(j_1,\ldots, j_k)\in\F_1\setminus\B} \pr\left(\cup_{i=1}^k H_{j_i}\subseteq \G(n,d)\right) \sim  \left(\frac{(n)_3}{6}\right)^k \lambda^k=\mu^k. \label{sum_F1_1} 
\end{equation}
On the other hand, by Corollary~\ref{cor:conditional}, for any $(j_1,\ldots, j_k)\in \F_1\cap \B$, we have
\begin{equation}
\ex \left(\prod_{i=1}^k X_{H_{j_i}}\right) \le \lambda^k \left(1+O\left(\frac{k}{n}+\frac{k^2}{dn}+\frac{d^2k}{n^2}\right)\right) \sim \lambda^k,\label{simple_kth_moments}
\end{equation}
and again by Claim~\ref{claim:F1} we have
\begin{equation}
\sum_{(j_1,\ldots, j_k)\in\F_1\cap\B} \pr\left(\cup_{i=1}^k H_{j_i}\subseteq \G(n,d)\right) = o\left(\left(\frac{(n)_3}{6}\right)^k \right)\lambda^k=o(\mu^k). \label{sum_F1_2}
\end{equation}
Now Lemma~\ref{lem:F1} follows by~\eqn{sum_F1_1} and~\eqn{sum_F1_2}.

Finally, we prove Claim~\ref{claim:F1}. First, we prove 
\begin{equation}
|\F_1\cap\B|=o\left(\left(\frac{(n)_3}{6}\right)^k\right).\label{B_bound}
\end{equation}

For each $j\ge 1$ let $\C_j(i)$ be the set of $k$-tuples with $h_j=i$. Then, $[N]^k=\cup_{i\ge 0} \C_j(i)$ for each $j$. Define a switching from $\C_j(i)$ to $\C_j(i-1)$ as follows. Choose $v$ and $i_1<i_2<\cdots<i_j$ such that  $h(v)=j$ and $H_{i_1},\ldots, H_{i_j}$ all use $v$. There are exactly $i$ ways to choose them. Then, we choose $j-1$ distinct vertices $v_1<\cdots<v_{j-1}$, such that $h(v_{\ell})=0$ for all $1\le \ell\le j-1$. Modify  $H_{i_{\ell}}$ for each $2\le \ell\le j$ so that $v$ is replaced by $v_{\ell-1}$. The resulting $k$-tuple is in $\C_j(i-1)$. The number of ways to perform a switching is at least $i(n-3k)_{j-1}$. On the other hand, the number of inverse switchings is at most $(k)_j$. Hence, as $k=O(d^{3/2})$,
\[
\frac{|\C_j(i)|}{|\C_j(i-1)|}\le \frac{(k)_j}{i(n-3k)_{j-1}} \le \frac{k^j}{i (n-3k)^{j-1}}=O\left(\frac{d^{3j/3}}{i n^{j-1}}\right).
\]
As $d=o(n^{1/2})$ it follows that
\[
\left|\cup_{j\ge 4} \cup_{i\ge 1}\C_j(i)\right| \le \sum_{j\ge 4} \sum_{i\ge 1} |\C_j(i)|= O(d^6/in^3) |[N]^k| =o(|[N]^k|).
\]
For $j=2$,
\[
\frac{|\C_2(i)|}{|\C_2(i-1)|}\le \frac{(k)_2}{i(n-3k)} =O\left( \frac{d^3}{i n}\right).
\]
It follows that the number of $k$-tuples where $h_2 \ge d^2/\sqrt{n}$ is $o(|[N]^k|)$ as $d=o(\sqrt{n})$. Similarly, the number of $k$-tuples where $h_3 \ge d^2/\sqrt{n}$ is $o(|[N]^k|)$. This confirms~\eqn{B_bound}.

Next, we prove $|\F_1|\sim |[N]^k|$ which together with~\eqn{B_bound} confirms Claim~\ref{claim:F1}. The total number of $k$-tuples in $[N]^k$ is $\left(\frac{(n)_3}{6}\right)^k$. It is thus sufficient to show that 
\begin{equation}
\Big|[N]^k\setminus \F_1\Big|=o\left(\left(\frac{(n)_3}{6}\right)^k\right). \label{F1_complement_bound}
\end{equation}
Let $\F_1'\subseteq [N]^k$ be the set of $(j_1,\ldots, j_k)$ where there exists $j_{u}<j_{v}$ such that $H_{j_u}$ and $H_{j_v}$ share at least an edge. Then,  $\F_1'=[N]^k\setminus \F_1$. Obviously,
\[
|\F_1'|\le k^2 n^4 \left(\frac{(n)_3}{6}\right)^{k-2} = \left(\frac{(n)_3}{6}\right)^k \cdot O\left(\frac{k^2}{n^2}\right)=o\left( \left(\frac{(n)_3}{6}\right)^k\right),
\]
as $k^2=O(d^3)=o(n^2)$. This completes the proof for Claim~\ref{claim:F1}. \qed\smallskip

{\em Proof of Lemma~\ref{lem:F2}.\ } Given $\bfj=(j_1,\ldots, j_k)\in\F_2$, let $\I_j(\bfj)$ be the set of $i$ where $H_{j_i}$ intersects exactly $j$ edges with $\cup_{\ell<i} H_{j_{\ell}}$, for $j\in\{1,2,3\}$. By the definition of $\F_2$, 
$|\I_1(\bfj)|+|\I_2(\bfj)|+|\I_3(\bfj)|\ge 1$ for any $\bfj\in\F_2$. 
We will consider a set of parameters for $\bfj$ and count elements in $\F_2$ according to these parameters. The first three parameters are
\[
i=|\I_1(\bfj)|,\quad j=|\I_2(\bfj)|,\quad \ell=|\I_3(\bfj)|.
\]
Let $u_i,v_i,w_i$ be the three vertices occupied by the triangle $H_{j_i}$. We say that $z_i$, $z\in\{u,v,w\}$, creates a collision, if $z_i$ has been occupied by some triangle $H_{j_c}$, $c<i$. $H_{j_i}$ can create at most 3 collisions. Suppose $i\in \I_1(\bfj)$ and $u_iv_i$ is an edge that has been occupied by some triangle $H_{j_c}$, $c<i$, then immediately $u_i$ and $v_i$ each creates a collision. These collisions are called {\em inherent}. Suppose $i\in \I_2(\bfj)$ then $H_{j_i}$ creates 3 collisions all of which are called {\em inherent}. All other collisions are called noninherent. 

Let
\begin{align*}
t\quad &\mbox{denote the number of noninherent collisions. }
\end{align*}
Finally, let
\begin{align*}
\F_2(i,j,\ell,t)&\ \mbox{be the set of $\bfj$ satisfying the set of parameters};\\
\F_2(i,j,\ell)&=\cup_{t<\log n\cdot k^2/n} \F_2(i,j,\ell,t).
\end{align*}
Now 
\[
\F_2=\left(\cup_{t\ge \log n\cdot k^2/n} \cup_{(i,j,\ell): i+j+\ell\ge 1} \F_2(i,j,\ell,t)\right)\cup\left(\cup_{(i,j,\ell): i+j+\ell\ge 1} \F_2(i,j,\ell)\right)
\]

Fix $\bfj\in \F_2$.  We show that
\begin{equation}
\pr\left(\cup_{i=1}^k H_{j_i}\subseteq \G(n,d)\right) \le \lambda^{k-i-j-\ell} \left(\frac{d}{n}\right)^{2i+j}\left(1+O\left(\frac{k}{n}+\frac{k^2}{dn}+\frac{d^2 k}{n^2}\right)\right).\label{Fij_bound}
\end{equation}
Since
\[
\pr\left(\cup_{i=1}^k H_{j_i}\subseteq \G(n,d)\right)=\prod_{i=1}^k \pr\left(H_{j_i}\subseteq \G(n,d) \mid (\cup_{\ell<i} H_{j_{\ell}})^+\right).
\]
By Corollary~\ref{cor:conditional}, besides the relative error $O(1/n+k/dn+d^2/n^2)$, if $i\notin (\I_1(\bfj)\cup \I_2(\bfj)\cup \I_3(\bfj))$ then the above conditional probability is at most $\lambda$; if $i\in \I_1(\bfj)$ then the conditional probability is at most $(d/n)^2$; if $i\in \I_2(\bfj)$ then the conditional probability is at most $d/n$; and if $i\in \I_3(\bfj)$ then the conditional probability is 1. Multiplying them together yields the main term in~\eqn{Fij_bound}. There are $O(k)$ terms in the product, and multiply the relative errors in each term gives the error as in~\eqn{Fij_bound}. As $k=O(d^{3/2})$ and $d=o(\sqrt{n})$ the error in~\eqn{Fij_bound} is $o(1)$.

Next, we prove that for some absolute constant $C>0$:
\begin{align}
|\F_2(i,j,\ell,t)|&\le \frac{k^{i+j+\ell}}{i! j!\ell!} (3kn)^i (6kd)^j X^{\ell}\binom{3k}{t}\left(\frac{(n)_3}{6}\right)^{k-i-j-\ell} \left(\frac{Ck}{n}\right)^{t}, \label{F2ij_bound}\\
|\F_2(i,j,\ell)|&\le \frac{k^{i+j+\ell}}{i! j!\ell!} (3kn)^i (6kd)^j Y^{\ell} \left(\frac{(n)_3}{6}\right)^{k-i-j-\ell}, \label{F2ij_bound2}
\end{align}
where
\[
X=C(i^{3/2}+j^{3/2}+t^{3/2}),\quad Y=C(i^{3/2}+j^{3/2}+(\log n\cdot k^2/n)^{3/2}).
\]
There are at most $k^{i+j+\ell}/i! j!\ell!$ ways to fix $\I_1$ for $\I_1(\bfj)$, $\I_2$ for $\I_2(\bfj)$, and $\I_3$ for $\I_3(\bfj)$. For every $h\in \I_1$, there are at most $3k$ ways to choose the edge in $H_{j_h}$ which intersect with $\cup_{c<h} H_{j_{c}}$, and then at most $n$ ways to choose the other vertex in $H_{j_h}$. For every $h\in\I_2$, there are at most $3k$ ways to choose an edge $x$ in $\cup_{c<h} H_{j_{c}}$, and then at most $2d$ ways to choose another edge incident with $x$ in $\cup_{c<h} H_{j_{c}}$. The upper bound $2d$ follows as we may assume that every vertex is incident with at most $d$ edges of $\cup_{i=1}^k H_{j_i}$ --- otherwise, $\pr(\cup_{i=1}^k H_{j_i} \subseteq \G(n,d))=0$. We claim the following bound on the number of holes.
\begin{claim}\label{c:holes}
The number of holes in $\F_2(i,j,\ell,t)$ is at most $X$.
\end{claim}
{\em Proof of Claim~\ref{c:holes}. } First we prove that any graph with $x$ edges has at most $3x^{3/2}$ triangles. Let $t_1,\ldots t_n$ be the degree sequence of the graph, and we call a vertex big if its degree is at least $\sqrt{x}$ and small otherwise. Then the number of big vertices is at most $2x/\sqrt{x}=2\sqrt{x}$. A triangle in the graph either uses only big vertices, or occupies a small vertex. Hence, the number of triangles in the graph is at most
\[
\binom{2\sqrt{x}}{3}+\sum_{i: t_i\le\sqrt{x}} \binom{t_i}{2} \le \frac{4}{3}x^{3/2} + \frac{2x}{\sqrt{x}} \frac{(\sqrt{x})_2}{2}\le 3x^{3/2},
\] 
where the maximum of the sum above appears when all non-zero degrees of small vertices are equal to $\sqrt{x}$.

Now, we bound the holes given parameters $i$, $j$ and $t$. Let $H=\cup_{h=1}^k H_{j_h}$ and we bound the largest possible number of holes that $H$ may contain. If $h\in \I_1\cup\I_2$ then colour all three edges in $H_{j_h}$ red. Then the number of red edges is at most $3(i+j)$. Let $T_1$ be the number of holes where all their three edges are red, and let $T_2$ be the number of other holes.
Immediately, $T_1\le 3(i+2j)^{3/2}\le C(i^{3/2}+j^{3/2})$ and $T_2\le Ct^{3/2}$. The claim follows. \qed 

First consider $\F_2(i,j,\ell,t)$ for $t\ge \log n\cdot k^2/n$.
By the above claim, there are at most $X^{\ell}$ ways to fix all $j_h$ for $h\in \I_3$.   Finally, there are at most $\left((n)_3/6\right)^{k-i-j-\ell}$ ways to fix $(j_i)_{i\notin (\I_1\cup \I_2\cup\I_3)}$. There are at most $\binom{3k}{t}$ ways to choose the set of noninherent collisions, and each noninherent collision contributes an additional $Ck/n$ factor because there are at most $3k$ instead of $n$ choices to choose the vertex where a noninherent collision occurs. Multiplying all these upper bounds together yields~\eqn{F2ij_bound}.

For~\eqn{F2ij_bound2}, the argument is similar except that we do not take union bound over all possible vertex collisions. Instead, we use the fact that $X\le Y=C(i^3+j^3+(\log n\cdot k^2/n)^{3/2})$ when $t<\log n\cdot k^2/n$. Thus the number of ways to choose $j_h$ where $h\in \I_3$ is at most $Y^{\ell}$ for all $\bfj\in\F_2(i,j,\ell)$, and~\eqn{F2ij_bound2} follows. 

By~\eqn{Fij_bound} and~\eqn{F2ij_bound} we have
\begin{eqnarray*}
&&\sum_{(j_1,\ldots, j_k)\in\F_2(i,j,\ell,t)} \pr\left(\cup_{i=1}^k H_{j_i}\subseteq \G(n,d)\right)\le  (1+o(1)) |\F_2(i,j,\ell,t)| \cdot  \lambda^{k-i-j-\ell} \left(\frac{d}{n}\right)^{2i+j}\\
&& \hspace{1cm} \le (1+o(1)) \frac{k^{i+j+\ell}}{i! j!\ell!} \left(\frac{(n)_3}{6}\right)^{k-i-j-\ell} (3kn)^i (6kd)^j X^{\ell}  \left(\frac{Ck^2}{tn}\right)^{t} \cdot\lambda^{k-i-j-\ell} \left(\frac{d}{n}\right)^{2i+j}\\
&& \hspace{1cm} =O( \mu^k)\cdot  \left(\frac{Ck^2d^2}{i\mu n}\right)^i \left(\frac{Ck^2d^{2}}{j\mu n}\right)^j \left(\frac{CkX}{\ell\mu }\right)^{\ell}\left(\frac{Ck^2}{tn}\right)^{t} \label{bound1}\\
&& \hspace{1cm} =O( \mu^k)\cdot    \left(\frac{Ck^2d^2}{i\mu n}\right)^i \left(\frac{Ck^2d^{2}}{j\mu n}\right)^j \exp\left(\frac{CkX}{\mu }\right) \left(\frac{Ck^2}{tn}\right)^{t}
\end{eqnarray*}
by by defining $(x/0)^0=1$, and noting that $(x/\ell)^{\ell} \le e^x$.
It follows then that
\begin{eqnarray}
&&\sum_{(j_1,\ldots, j_k)\in\F_2(i,j,\ell,t)} \pr\left(\cup_{i=1}^k H_{j_i}\subseteq \G(n,d)\right)\nonumber\\
&&\hspace{1cm} =O( \mu^k) \cdot  \left(\frac{Ck^2d^2}{i\mu n}\right)^i \left(\frac{Ck^2d^{2}}{j\mu n}\right)^j   \exp\left(\frac{Ck(i^{3/2}+j^{3/2}+b^3+t^{3/2})}{\mu }\right)\left(\frac{Ck^2}{tn}\right)^{t}\nonumber\\
&&\hspace{1cm} =O( \mu^k) \cdot \exp\left(i\left(\ln(Ck^2d^2/i\mu n)+\frac{Ck\sqrt{i}}{\mu}\right)+j\left(\ln(Ck^2d^2/j\mu n)+\frac{Ck\sqrt{j}}{\mu}\right) \right)\left(\frac{Ck^2}{tn}\right)^{t}.\nonumber\\
\label{sum0}
\end{eqnarray}
Since $\mu=\Theta(d^3)$, $k=\Theta(d^{3/2})$, $d=o(n^{1/2})$, and $t\ge \log n\cdot k^2/n$, it follows immediately that 
\begin{equation}
\sum_{t\ge \log n\cdot k^2/n}\sum_{i,j,\ell}\sum_{(j_1,\ldots, j_k)\in\F_2(i,j,\ell,t)} \pr\left(\cup_{i=1}^k H_{j_i}\subseteq \G(n,d)\right)=o(\mu^k).\label{sum1}
\end{equation}

Now combining~\eqn{Fij_bound} and~\eqn{F2ij_bound2} we have
\begin{eqnarray}
&&\sum_{(j_1,\ldots, j_k)\in\F_2(i,j,\ell)} \pr\left(\cup_{i=1}^k H_{j_i}\subseteq \G(n,d)\right)\le  (1+o(1)) |\F_2(i,j,\ell)| \cdot  \lambda^{k-i-j-\ell} \left(\frac{d}{n}\right)^{2i+j}\nonumber \\
&& \hspace{1cm} \le (1+o(1)) \frac{k^{i+j+\ell}}{i! j!\ell!} \left(\frac{(n)_3}{6}\right)^{k-i-j-\ell} (3kn)^i (6kd)^j (C(i^{3/2}+j^{3/2}+(\log n\cdot k^2/n)^{3/2}))^{\ell} \nonumber\\
&& \hspace{1.5cm}   \cdot\lambda^{k-i-j-\ell} \left(\frac{d}{n}\right)^{2i+j}\nonumber\\
&& \hspace{1cm} =O( \mu^k)\cdot  \left(\frac{Ck^2d^2}{i\mu n}\right)^i \left(\frac{Ck^2d^{2}}{j\mu n}\right)^j \left(\frac{Ck(i^{3/2}+j^{3/2})}{\ell\mu }+C\left(\frac{d\log n}{n}\right)^{3/2}\right)^{\ell}\label{sum2}
\end{eqnarray}
Since $\mu=\Theta(d^3)$, $k=\Theta(d^{3/2})$, $d=o(n^{1/2})$, summing the above over all $(i,j,\ell)$ where $i+j+\ell\ge 1$ yields $o(\mu^k)$. \qed

\subsection{Distribution of $Z_{C_3}$ for $d=\Theta(\sqrt{n})$}
\label{sec:large}

Again let $N=\binom{n}{3}$ and let $H_1,\ldots, H_N$ be an enumeration of all copies of triangles in $K_n$. Let $\F$ be the set of  ordered $k$-tuples $(j_1, \ldots, j_k)\in [N]^k$ that are pairwise distinct. We want to estimate the following sum for $k=O(d^{3/2})$.
\begin{equation}
\ex (Z)_k = \sum_{(j_1,\ldots, j_k)\in\F} \pr\left(\cup_{i=1}^k H_{j_i}\subseteq \G(n,d)\right). \label{factorial}
\end{equation}

We partition the summation into parts $\F_1$ and $\F_2=\F\setminus \F_1$, where $\F_1\subseteq \F$ is the set of $\bfj$ satisfying one the following:
\begin{itemize}
\item[(a)]  $|\I_3(\bfj)|\ge 1$, or $h(v)> \log^2 n$ for some $v$.
\item[(b)]  there  exists $H_{j_i}$ which intersects at least two edges with $\cup_{\ell<i} H_{j_{\ell}}$.
\item[(c)]  there are more than $\log n$ $i$'s  such that $H_{j_i}$  intersects exactly one edge with $\cup_{\ell<i} H_{j_{\ell}}$.
\item[(d)]  there exists an edge contained in at least three triangles.
\item[(e)] there is some $H_{j_i}$ such that at least two of its edges are each contained in some other triangles.
\end{itemize}

We will prove the following.
\begin{lemma}\label{lem2:F1}
\[
\sum_{\bfj\in\F_1} \pr\left(\cup_{i=1}^k H_{j_i}\subseteq \G(n,d)\right) =o(\mu^k).
\]
\end{lemma}

\begin{lemma}\label{lem2:F2}
\[
\sum_{\bfj\in\F_2} \pr\left(\cup_{i=1}^k H_{j_i}\subseteq \G(n,d)\right) \sim \mu^k.
\]
\end{lemma}


{\em Proof of Lemma~\ref{lem2:F2}.\ } The structure of the union of the $k$ triangles $H_{j_1},\ldots, H_{j_k}$ are simple, for $\bfj\in\F_2$. Except for at most $\log n$ pairs of triangles, each pair sharing exactly one edge, all the other triangles are edge disjoint. Moreover, every vertex is contained in at most $\log^2 n$ triangles.  However, as we will show, there are $\Omega(d)$ vertices $v$ with $h(v)\ge 2$. This means that the joint probability of a set of $k$ edge disjoint triangles will not be asymptotic to $\lambda^k$ any more, which was the case in the proof of Lemma~\ref{lem:F1}.

We will partition $\F_2$ according to two parameters. An edge is called a {\em double} edge, if it is contained in two triangles. For ${\bfj}\in \F_2$, the number of double edges is exactly $|\I_1(\bfj)|$. Recall that $h_i$ denote the number of vertices $v$ with $h(v)=i$. Modify ${\cal B}$ to state the set of $\bfj$ where
\[
|h_2-9k^2/2n|> d^{2/3},\ \mbox{or}\ h_3+h_4> d \cdot d^{7/2}\log n/n^2,\ \mbox{or} \ h_j>0\ \mbox{for some $j\ge 5$.}
\]

 Define
\begin{eqnarray*}
\F_2(i)=\{\bfj\in \F_2:\ |\I_1(\bfj)|=i\},&& \F_2^0(i)=\F_2(i)\setminus \B, \quad \F_2^1(i)=\F_2(i)\cap \B \\
\F_2^0=\cup_{i=0}^{\log n} \F_2^0(i),&& \F_2^1=\cup_{i=0}^{\log n} \F_2^1(i)=\F_2\setminus \F_2^0.
\end{eqnarray*}
Then for every $i\le \log n$,
\begin{eqnarray*}
|\F_2(i)| &=&\binom{k}{2i} \frac{(2i)!}{2^i i!} \left(\frac{(n)_3}{6}-O(kn)\right)^{k-2i} \left(\frac{(n-O(i))^4}{2}\right)^i \\
&=& \frac{k^{2i}}{4^i i!} \left(\frac{(n)_3}{6}\right)^{k-2i}n^{4i}\exp\left(O\left(\frac{k^2}{n^2}+\frac{i^2}{k}\right)\right)\\
&\sim & \frac{k^{2i} n^{4i}}{4^i i!} \left(\frac{(n)_3}{6}\right)^{k-2i}.
\end{eqnarray*}
\begin{claim}\label{claim:h2}
\[
|\F_2^0(i)|\sim |\F_2(i)|\ \mbox{for all $i\le \log n$, and}\quad \sum_{\bfj\in\F_2^1}   \pr\left(\cup_{i=1}^k H_{j_i}\subseteq \G(n,d)\right)= o(\mu^k).
\]
\end{claim}

Next, fix $\bfj\in \F_2^0(i)$.  We show that
\begin{eqnarray}
\pr\left(\cup_{i=1}^k H_{j_i}\subseteq \G(n,d)\right) &\sim & \lambda^{k-2i}  \left(\frac{d}{n}\right)^{5i} \exp\left( -\frac{9k^2}{dn}\right).\label{triangle_prob_error}
\end{eqnarray}
The proof is almost the same as for~\eqn{kthMoments}. The accumulative error, except for that from $\phi_F(uv)/dn$, is $O\left((d^2 k+k^2)/n^2+k^3/d^2n^2\right)=o(1)$. The only two differences are listed below:
\begin{itemize}
\item[(a)] The term $-2|F|$ in $\phi_F(uv)$ in~\eqn{phi2} is no longer negligible. Its accumulative contribution after taking the product of the conditional probabilities is
\[
\prod_{\ell=1}^{3k-O(i)} \left(1+\frac{2\ell}{dn}\right) = \exp\left( \frac{9k^2}{dn} +O\left(\frac{ik}{dn}+\frac{k^3}{d^2n^2}\right)\right)\sim \exp\left( \frac{9k^2}{dn}\right).
\] 
\item[(b)] For each $v$ where $h(v)=2$, the joint probability of the two triangles hitting $v$ cannot be approximated by $\lambda^2$ with negligible error because $h_2=\Theta(d)$. Their accumulative contribution  after taking the product of the conditional probabilities becomes
\[
\left(\frac{(d-2)(d-3)}{d(d-1)}\right)^{h_2}=\exp\left(\left(-\frac{4}{d}+O(1/d^2)\right)\left(\frac{9k^2}{2n}+O(d^{2/3})\right)\right)\sim \exp\left(-\frac{18k^2}{dn}\right).
\]
\end{itemize}
Multiplying the additional factors in (a) and (b) yields~\eqn{triangle_prob_error}.
Finally,
\begin{eqnarray}
\sum_{\bfj\in\F_2^0} \pr\left(\cup_{i=1}^k H_{j_i}\subseteq \G(n,d)\right)&=&\sum_{i=0}^{\log n} \sum_{\bfj\in\F_2^0(i)} \pr\left(\cup_{i=1}^k H_{j_i}\subseteq \G(n,d)\right)\nonumber\\
&\sim&  \sum_{i=0}^{\log n} \left(\frac{(n)_3}{6}\right)^{k-2i} \frac{k^{2i}n^{4i}}{4^ii!} \lambda^{k-2i}  \left(\frac{d}{n}\right)^{5i} \exp\left(- \frac{9k^2}{dn}\right)\nonumber\\
&\sim& \mu^k \exp\left(-\frac{9k^2}{dn}\right) \sum_{i=0}^{\log n} \frac{6^{2i}}{n^{6i}} \frac{k^{2i} n^{4i}}{4^ii!}\left(\frac{n}{d}\right)^i\nonumber\\
&\sim& \mu^k \exp\left(-\frac{9k^2}{dn}\right) \exp\left(\frac{9k^2}{dn}\right)\nonumber\\
&=& \mu^k, \label{F2_0}
\end{eqnarray}
as $k=O(d^{3/2})$, $d=\Theta(\sqrt{n})$ and $\mu=\Theta(d^3)$. \smallskip

It only remains to prove Claim~\ref{claim:h2}.

{\em Proof of Claim~\ref{claim:h2}.\ } The proof is basically the same as the proof for Claim~\ref{claim:F1}, except that in the switching we avoid switch away the double edges --- there are at most $\log n$ of them, and we estimate the number of switchings more carefully for $h_2$. The argument for bounds on $h_j$, $j\ge 3$ is much simpler and we skip the details. Now consider $h_2$, and assume that $h_3+h_4=o(d)$ and $h_j=0$ for all $j\ge 5$. Let $\C(i,j)$ be the set of $\bfj\in \F_2$ with $i$ double edges with $h_2=j$. Let $\C(i)=\cup_{j\ge 0} \C(i,j)$.
Define a switching from $\C(i,j)$ to $\C(i,j-1)$ as follows. Choose $v$ and $i_1<i_2$ such that  $h(v)=2$, $v$ is not incident with a double edge, and $H_{i_1} $ and $H_{i_2}$ both hit $v$. There are exactly $j-O(i)$ ways to choose them. Then, we choose a vertex $u$ such that $h(u)=0$. Modify  $H_{i_{2}}$  so that $v$ is replaced by $u$. The resulting $k$-tuple is in $\C(i,j-1)$. The number of ways to perform a switching is  $(j-O(i))(n-3k+O(\sum_{j} h_j))$. On the other hand, the number of inverse switchings is $9\binom{k-2(j-O(i))-\sum_{j\ge 3} j h_j}{2} $. Hence, 
\begin{eqnarray}
\frac{|\C(i,j)|}{|\C(i,j-1)|}&=& \frac{9\binom{k-2(j-O(i))-\sum_{j\ge 3} j h_j}{2}}{j(n-3k+O(\sum_{j} h_j))}=\frac{9(k-2j)^2}{2j(n-3k)}\left(1+O\left(\frac{1}{\sqrt{dn}}+\frac{d}{n}+\frac{i}{j}\right)\right)\nonumber\\
&=&\frac{9(k-2j)^2}{2j(n-3k)}\left(1+O\left(\frac{\log n}{\sqrt{n}}\right)\right).\label{ratio_h2}
\end{eqnarray}
Let $j^*$ be the root of $9(k-2j)^2=2(n-3k)$. Then,
\[
j^*=\frac{5k}{12}+\frac{n}{36}-\frac{\sqrt{n^2+30kn-99k^2}}{36}=\frac{9k^2}{2n}+O\left(\sqrt{d^3/n}\right).
\]
Now using~\eqn{ratio_h2} and following standard calculations we have
$
\sum_{|j-j^*|>d^{2/3}} |\C(i,j)| =o(|\C(i)|)$ and the first assertion of the claim follows. The second assertions follows by the same calculations as for~\eqn{triangle_prob_error} and~\eqn{F2_0}, except that we may use the bound in Corollary~\ref{cor:conditional} which is of simpler form.
For each $\bfj\in \F_2\cap\B$, 
\[
\pr\left(\cup_{i=1}^k H_{j_i}\subseteq \G(n,d)\right) \le \lambda^{k-2i}  \left(\frac{d}{n}\right)^{5i} \left(1+O\left(\frac{k}{n}+\frac{k^2}{dn}+\frac{d^2k}{n^2}\right)\right)=O\left(\lambda^{k-2i}  \left(\frac{d}{n}\right)^{5i}\right).
\]
Since $|\F^1_2(i)| = o(|\F_2(i)|$ for all $i\le\log n$,
\begin{eqnarray*}
\sum_{\bfj\in\F_2^1} \pr\left(\cup_{i=1}^k H_{j_i}\subseteq \G(n,d)\right)&=&\sum_{i=0}^{\log n} o(|\F_2(i)|) \cdot O\left(\lambda^{k-2i}  \left(\frac{d}{n}\right)^{5i}\right)\\
&=&o(1) \sum_{i=0}^{\log n} \frac{k^{2i} n^{4i}}{4^i i!} \left(\frac{(n)_3}{6}\right)^{k-2i} \lambda^{k-2i}  \left(\frac{d}{n}\right)^{5i},
\end{eqnarray*}
which is $o(\mu^k)$ by following the same calculations as in~\eqn{F2_0}. \qed
 \smallskip

{\em Proof of Lemma~\ref{lem2:F1}.\ } 
 Let $\I_j(\bfj)$ be as defined before. Let $\F_1^1$ be the set of $\bfj\in\F_1$ where $\bfj$ induces a hole.  With almost the same proof as in Lemma~\ref{lem:F2}, we consider $t$: the number of noninherent collisions. For $t\ge \log n\cdot k^2/n$, we get~\eqn{sum0}, and thus~\eqn{sum1}. Note that the change of the range of $d$ from $o(\sqrt{n})$ to $O(\sqrt{n})$ does not affect the analysis, as $t\ge \log n\cdot k^2/n$ and thus the term $(Ck^2/tn)^t$ dominates the rate the whole product vanishes. 
 For $t<\log n\cdot k^2/n$, we get~\eqn{sum2}. Summing over all $(i,j,\ell)$
 where $\ell\ge 1$ yields $o(\mu^k)$.


Next, let $\F_1^2\subseteq \F_1\setminus \F_1^1$ be the set of $\bfj$ where $i+j>\log n$ where $i=|\I_1(\bfj)|$ and $j=|\I_2(\bfj)|$. As $\bfj\notin\F_1^1$ we know $\I_3(\bfj)$ must be empty.
 Let 
 \[
 \F_1^2(i,j)=\{{\bfj}:\ \I_1(\bfj)=i,\ \I_2(\bfj)=j \}.
 \]
We use the bound~\eqn{F2ij_bound2} (with $\ell=0$) for $|\F_1^2(i,j)|$:
\[
|\F_1^2(i,j)|\le \frac{k^{i+j}}{i! j!} \left(\frac{(n)_3}{6}\right)^{k-i-j} (3kn)^i (6kd)^j.
\]
Thus,
\begin{eqnarray*}
&&\sum_{i,j: i+j>\log n}\sum_{\bfj\in\F_1^2(i,j)} \pr\left(\cup_{i=1}^k H_{j_i}\subseteq \G(n,d)\right)\le  O(1)\cdot \mu^k \sum_{i,j: i+j\ge \log n} \frac{(3k^2d^2/\mu n)^i}{i!} \frac{(6k^2d^{2}/\mu n)^j}{j!}=o(\mu^k),
\end{eqnarray*}
as $k=O(d^{3/2})$ and thus $3k^2d^2/\mu n=O(1)$.

 Next, let $\F_1^3\subseteq \F_1\setminus (\F_1^1\cup\F_1^2)$ be the set of $\bfj$ where $h(v)>\log^2 n$. The summation over $\F_1^3$ can be treated similarly as in $\F_1^1$.  
 Fix $v\in [n]$. Since $\bfj$ induces no holes, $\I_3(\bfj)=0$. Let ${\cal L}(\bfj)$ be the set of $i$ such that $H_{j_i}$ hits $v$. 
Let 
\begin{eqnarray*}
\F_1^3(i,j,r)&=&\{{\bfj}\in \F_1^3:\ \I_1(\bfj)=i,\ \I_2(\bfj)=j,\ \ {\cal L}(\bfj)=r\}.
\end{eqnarray*}

Since $\bfj\notin \F_1^2$, $i\le \log n$ and $j\le \log n$. Let $r>\log^2 n$.  There are at most
$k^{i+j+r}/i!j!r!$ ways to fix $\I_1$, $\I_2$ and ${\cal L}$ for $\I_1(\bfj)$, $\I_2(\bfj)$ and ${\cal L}(\bfj)$ respectively. Given $\I_1$, $\I_2$ and ${\cal L}$, let ${\cal L}'={\cal L}\setminus (\I_1\cup \I_2)$ and we must have $|{\cal L}'|\ge r-2\log n=r(1-o(1))$ as $r>\log^2 n$. Thus,
\[
|\F(i,j,r)|\le  \frac{k^{i+j+r}}{i! j!r!} \left(\frac{(n)_3}{6}\right)^{k-i-j} (3kn)^i (6kd)^j n^{2r},
\]
as for each $r\in {\cal L}'$, there are at most $n^2$ ways to choose the two vertices other than $v$. Hence,
\begin{eqnarray*}
&&\sum_{\substack{i,j,r:\\ r>\log^2 n}}\sum_{\bfj\in\F_1^3(i,j,\ell)} \pr\left(\cup_{i=1}^k H_{j_i}\subseteq \G(n,d)\right)\\
&&\hspace{1cm}= O(1)\cdot \sum_{i,j,r: r>\log^2 n} \frac{k^{i+j+r}}{i! j!r!} \left(\frac{(n)_3}{6}\right)^{k-i-j-r(1-o(1))} (3kn)^i (6kd)^j n^{2r} \lambda^{k-i-j}\left(\frac{d}{n}\right)^{2i+j}\\
&&\hspace{1cm}= O(\mu^k) \sum_{r>\log^2 n} \frac{(k/n^{1-o(1)})^r}{r!} \sum_{i\ge 0} \frac{(3k^2d^2/\mu n)^i}{i!} \sum_{j\ge 0}\frac{(6k^2d^{2}/\mu n)^j}{j!} =o(\mu^k),
\end{eqnarray*}
as 
$k^2d^2/\mu n = O(1)$ 
and 
$k/n^{1-o(1)}=o(1)$.

Next, let $\F_1^4\subset \F_1\setminus (\F_1^1\cup\F_1^2\cup\F_1^3)$ be the set of ${\bfj}$ where $\I_2(\bfj)$ is nonempty.
Let 
\[
\F_1^4(i,j) = \{\bfj\subseteq \F_1^4:\ \I_1(\bfj)=i,\ \I_2(\bfj)=j\}.
\] 
As now $\bfj\notin \F_1^3$, $h(v)\le \log^2 n$ for every $n$ and thus we have the following better upper bound than~\eqn{F2ij_bound} for $|\F_1^4(i,j)|$:
\[
|\F_1^4(i,j)| \le \left(\frac{(n)_3}{6}\right)^{k-i-j} (3kn)^i (12k\log^2 n)^j.
\]
This is because every vertex in $\cup_{i=1}^{k} H_{j_i}$ is incident with at most $2\log^2 n$ instead of $d$ edges.  
Using this improved upper bound, we have
\begin{eqnarray*}
&&\sum_{\bfj\in\F_1^4}\pr\left(\cup_{i=1}^k H_{j_i}\subseteq \G(n,d)\right) =O(\mu^k) \sum_{i\le \log n} \frac{(3k^2d^2/\mu n)^i}{i!} \sum_{j\ge 1} \frac{(6k^2d\log^2 n/\mu n)^j}{j!}=o(\mu^k).
\end{eqnarray*}

Finally, let $\F_1^5=\F_1\setminus(\cup_{i=1}^4\F_1^i)$, the set of $\bfj$ where either an edge is contained in more than 2 triangles, or there is a triangle with two edges each contained in some other triangle. 
Summing over $\bfj$ in the first case we have
\[
 \sum_{2\le i\le \log n}  \frac{k^{i}}{i!} \left(\frac{(n)_3}{6}\right)^{k-i} i \cdot i n \cdot (3kn)^{i-1}  \lambda^{k-i}\left(\frac{d}{n}\right)^{2i}=O(\mu^k \cdot \log^2 n/k)=o(\mu^k). 
\]
Now we consider the latter case. If $H_{j_u}$ has two edges each contained in some other triangles --- call them $H_{j_{\ell}}$ and $H_{j_r}$. Obviously $\ell\neq r$ since otherwise $H_{j_u}$ would be the same as $H_{j_{\ell}}$, contradicting with them being distinct. There are two sub-cases: (a) $\ell<i$ and $r>i$; (b) $r>\ell>i$. In either case, there are at most $i^3$ ways to choose $u,\ell$ and $r$. Then, similar to the above,
summing over $\bfj$ in the second case yields
\[
 \sum_{2\le i\le \log n}  \frac{k^{i}}{i!} \left(\frac{(n)_3}{6}\right)^{k-i} i^3\cdot O( n^2) \cdot (3kn)^{i-2}  \lambda^{k-i}\left(\frac{d}{n}\right)^{2i}=O(\mu^k \cdot \log^3 n/k)=o(\mu^k). 
\]
Combining the two cases we have $\sum_{\bfj\in\F_1^5}\pr\left(\cup_{i=1}^k H_{j_i}\subseteq \G(n,d)\right) =o(\mu^k)$, completing the proof of the lemma. \qed


\end{document}